\documentclass[11pt]{article}
\usepackage{graphics}
\usepackage{color}

\newfont{\frak}{eufm10 scaled\magstep1}
\newfont{\sfrak}{eufm8 scaled\magstep1}
\newfont{\bbb}{msbm10 scaled\magstephalf}
\newfont{\sbbb}{msbm7 scaled\magstephalf}

\def\D{\Delta}
\def\C{\mbox{\bbb{C}}}
\def\R{\mbox{\bbb{R}}}
\def\Z{\mbox{\bbb{Z}}}
\def\K{\mbox{\bbb{K}}}
\def\Q{\mbox{\bbb{Q}}}
\def\cd{\C^d}
\def\rd{\R^d}

\def\vz{\underline{z}}
\def\vx{\underline{x}}

\def\vu{\underline{u}}
\def\vw{\underline{w}}
\def\vrho{\underline{\rho}}
\def\vtheta{\underline{\theta}}
\def\ed{e_1,\ldots,e_d}
\def\ld{\lambda_1,\ldots,\lambda_d}
\def\xd{X_1,\ldots,X_d}
\def\d{\mbox{\frak d}}
\def\n{\mbox{\frak n}}
\def\ddu{\d^*}

\def\zd{(z_1,\cdots,z_d)}
\def\lorw{\longrightarrow}
\def\SC{\mbox{\sbbb{C}}}
\def\Dc{D_{\SC}}

\def\zset{\Psi^{-1}(0)}

\def\cono{C}

\def\regpsi{\Psi^{-1}(0)_{\D}}

\def\ut{\tilde{U}}
\def\uti{\tilde{U}_I}
\def\utj{\tilde{U}_J}
\def\ush{U^{\sharp}}
\def\uci{\widehat{U}_I}
\def\ucj{\widehat{U}_J}
\def\b{\mbox{\frak b}}

\def\ua{U_{\alpha}}
\def\uta{\tilde{U}_{\alpha}}
\def\fia{\phi_{\alpha}}

\def\vta{\tilde{V}_{\alpha}}
\def\wa{W_{\alpha}}
\def\wta{\tilde{W}_{\alpha}}
\def\pia{\psi_{\alpha}}
\def\ga{\Gamma_{\alpha}}
\def\ub{U_{\beta}}
\def\utb{\tilde{U}_{\beta}}

\def\fib{\phi_{\beta}}
\def\ga{\Gamma_{\alpha}}
\def\gb{\Gamma_{\beta}}
\def\da{\Delta_{\alpha}}
\def\gab{g_{\alpha\beta}}
\def\G{\Gamma}
\def\fa{f_{\alpha}}

\def\b{{\mbox{\frak{s}}}}
\def\sigmaf{\Sigma^{\diamond}_F}

\def\GF{G_{F}}

\def\squareforqed{\hbox{\rlap{$\sqcap$}$\sqcup$}}
\def\qed{\ifmmode\else\unskip\quad\fi\squareforqed}
\def\smartqed{\def\qed{\ifmmode\squareforqed\else{\unskip\nobreak\hfil
\penalty50\hskip1em\null\nobreak\hfil\squareforqed
\parfillskip=0pt\finalhyphendemerits=0\endgraf}\fi}}

\newtheorem{thm}{Theorem}[section]
\newtheorem{prop}[thm]{Proposition}
\newtheorem{lemma}[thm]{Lemma}
\newtheorem{defn}[thm]{Definition}
\newtheorem{remark}[thm]{Remark}
\newtheorem{example}[thm]{Example}
\newcommand{\proof}{\mbox{\textbf{ Proof.\ \ }}}

\title{\sc Convex polytopes and quasilattices from the symplectic viewpoint}
\author{\sc Fiammetta Battaglia\thanks{Partially supported by
GNSAGA (CNR).}
}
\date{Preliminary version, 15/04/03}
\date{}
\begin{document}
\maketitle
\begin{abstract}
We construct, for each convex
polytope, possibly nonrational and nonsimple, a family of compact
spaces that are stratified by quasifolds, i.e. 
each of these spaces is a collection of quasifolds glued together in an 
suitable way. A quasifold is a space
locally modelled on $\R^k$ modulo the action of a discrete,
possibly infinite, group. The way strata are glued to each other also
involves the action of an (infinite) discrete group.
Each stratified space is endowed with a
symplectic structure and
a moment mapping having the property that
its image gives the original polytope back.
These spaces may be
viewed as a natural generalization of symplectic toric varieties
to the nonrational setting.
\end{abstract}

\medskip

{\small 2000 \textit{Mathematics Subject Classification.} Primary: 53D05.
Secondary: 53D20, 32S60, 52B20.}

\medskip


\section*{Introduction}
Consider a vector space $\d$ of dimension $n$.
In \cite{delzant} Delzant shows that 
each $n$-dimensional simple 
convex polytope $\D\subset\d^*$, 
rational with respect to a lattice in $\d$ and satisfying
an additional integrality condition, is the image of a 
symplectic toric manifold via moment mapping. 
The result is obtained by providing
an explicit construction of such a manifold, 
which coincides in fact with 
the toric variety associated to the polytope and is naturally endowed with a 
symplectic structure whose symplectic volume is equal to the Euclidean 
volume of the polytope.
In \cite{p} Prato  extends Delzant's result 
to the nonrational setting, constructing, 
for each $n$-dimensional simple convex polytope $\D\subset\d^*$, 
{\it rational or not},
a family of $2n$-dimensional compact symplectic spaces,  
called {\it quasifolds}, 
naturally associated to $\D$.
Each space of the family admits the effective Hamiltonian action of
an $n$-dimensional {\it quasitorus} $D$, with a moment mapping
whose image is exactly $\D$.
A quasifold is a topological space locally modelled
on $\R^k$ modulo the action of a discrete, possibly infinite, group, and 
therefore it is not necessarily Hausdorff. 
A quasitorus is the natural analogue of a torus in this setting, it is
the group and quasifold $\d/Q$, where $Q$ is a quasilattice\footnote{ 
Quasilattices are quasiperiodic structures underlying quasiperiodic 
tilings and quasicrystals atomic order. There is now an extremely 
rich literature on the subject, whose origin goes back to the eighties for
quasicrystals
\cite{quasicristalli}--for an updated account see for example 
\cite{quasicristallirecente}, while aperiodic tilings were produced first in 
the sixties by R. Berger, followed by the works by Robinson \cite{rob}
and Penrose \cite{pen}.}.

The idea underlying the present paper is that in order  
to extend these results to arbitrary convex polytopes, neither simple, nor 
rational, quasifolds are still the natural structures but
a further degree of singularity has to be allowed, 
just as in the case of classical toric varieties, in which nonsimplicity
of the polytope brings in singularities which are not finite group 
quotient type.
We therefore introduce
spaces stratified by quasifolds: nonsimplicity of the polytope causes
the decomposition in strata of the corresponding topological spaces, whilst 
nonrationality of $\D$ produces the quasifold structure of the strata
and intervenes in the way strata are glued to each other, leading to
a definition of stratification that naturally extends the usual one.
More precisely 
we prove that to each
$n$-dimensional 
convex polytope $\D\subset\d^*$, there corresponds a family of compact
spaces that are stratified by symplectic quasifolds. 
Different members of the family correspond, like in \cite{p}, to
different choices of additional data attached to $\D$. 
Each space $M$ of the family admits the continuous effective action of an
$n$-dimensional quasitorus $D$ and a continuous mapping
$\Phi\,\colon\,M\lorw\d^*$ such that $\Phi(M)=\D$. The restriction
of the $D$-action to each stratum is smooth and Hamiltonian, with
moment mapping given by the restriction of $\Phi$.
The proof is based on the explicit construction of $M$
as symplectic quotient; the stratification
of $M$ is the one naturally induced by the decomposition by isotropy type, 
it is, in other words, the decomposition by singularity type and mirrors, 
as we shall see, the decomposition of the polytope in singular and 
nonsingular faces. Moreover, as in the rational case, the local structure 
of the stratification reflects the polytope shape. 
If we restrict to the rational case we obtain, from the point of 
view of singularities,
what expected from the classical theory of toric varieties, in addition
we gain further insight on the symplectic 
geometry of them; it should also be noticed that these spaces  
provide a wide range of explicit examples of symplectic stratified
spaces as defined in \cite{ls} (further details in 
Remark~\ref{choices}).

In the present paper, whose results have been announced, jointly with Elisa 
Prato, in \cite{bp2}, the spaces corresponding to arbitrary convex polytopes
are described in the symplectic setting.
They have a complex counterpart,
strictly related to the theory of toric varieties (for further details cf.
Remark~\ref{complex}; for the extension of the notion of toric variety
to simple nonrational convex polytope see \cite{cx}).  
This aspect will be treated in a subsequent paper \cite{nscx}.

Recent works by Karu and Bressler-Lunts deal with 
arbitrary convex polytopes from  
another viewpoint.
In order to
prove the Hard-Lefschetz theorem (\cite{karu}) and Hodge-Riemann bilinear
relations (\cite{bl}) for nonrational convex polytopes, they
make use of
a combinatorial construction, from the polytope data, of the intersection 
cohomology, thus bypassing the problem of the existence 
of a ``variety'' associated to $\D$. 
What we provide in here is on the contrary
an explicit construction of a geometric space, indeed of a whole
family of geometric spaces, that naturally correspond to $\D$.
The geometry and topology of our spaces and the relationship with 
properties of the corresponding polytope--volume, 
counting integer lattice points,
Hard Lefschetz--are all very natural questions related to our work.
A first step towards a better understanding of these different aspects, 
that will be pursued in the sequel, is to investigate 
cohomological properties  of our spaces.

The paper is organized as follows.
After a preliminary section in which we recall the basics about quasifolds,
we give in section~$2$ the definition of space stratified by quasifolds.
In section~$3$
we explicitly construct, adopting the procedure introduced by Delzant in 
\cite{delzant} and extended to the nonrational case by Prato in \cite{p},
a family of topological spaces naturally associated to a convex polytope $\D$.
In the last two  sections we describe the structure of these spaces, by
proving that they are spaces stratified 
by symplectic quasifolds. To guide the exposition we provide two 
model examples worked out in detail.

\vspace{.3cm}

{\bf Acknowledgments}. This work began in collaboration with Elisa Prato
and the results presented here have been announced in the joint paper
\cite{bp2}. The project developed in these works was initiated by
Prato's article \cite{p} and then carried on jointly in \cite{cx}. 
I am very grateful to Elisa Prato for having introduced me to the 
beautiful subject of quasifolds -- working together has been an enrichening 
experience.

\section{Preliminaries about quasifolds}

In this section we recall basic definitions about quasifolds,
for further details and for related notions, such as
symplectic forms and Hamiltonian actions, 
we refer the reader to Prato's article \cite{p}.

The local model for quasifolds
is a  manifold acted on diffeomorphically by a discrete group.
\begin{defn}[Model]\label{model}{\rm
    Let $\ut$ be a connected, simply connected manifold of
    dimension $k$ and let $\G$ be a discrete group acting smoothly
on $\ut$
    so that the set of points, $\ut_0$, where the action is free,
    is connected and dense. Consider the space of orbits, $\ut/\G$, 
of the action of the group
    $\G$ on the manifold $\ut$, endowed with the quotient topology, and the
    canonical projection $p\;\colon\;\ut\rightarrow \ut/\G$.
    A {\em model} of dimension $k$ is the triple $(\ut/\G,p,\ut)$, shortly
    $\ut/\G$.}
\end{defn}
\begin{remark}\label{discreto}{\rm
Nonclosed subgroups of Lie groups, for instance of tori,
play a central role in our construction.
Such groups are immersed Lie subgroups 
but of course they are not embedded. 
Most discrete groups
we are dealing with are of this kind, more precisely they are 
finitely generated
nonclosed subgroups of Lie groups--they are not of course discrete with
the induced topology.
}\end{remark}
\begin{remark}\label{simplyc}{\rm 
Let $(\ut/\G,p,\ut)$ be a triple with
$\ut$ not necessarily simply connected and satisfying all other requirements
of Definition~\ref{model}. 
We can obtain a model from $(\ut/\G,p,\ut)$ 
by the following procedure: 
consider the universal cover, $\pi\,\colon\,\ush\rightarrow\ut$, 
and its
fundamental group, $\Pi$. The manifold $\ush$ is connected and simply connected, the
mapping $\pi$ is smooth, the discrete group $\Pi$ acts smoothly, freely
and properly on the manifold $\ush$ and $\ut=\ush/\Pi$. 
Consider the extension of the
group $\G$ by the group $\Pi$, $1\longrightarrow \Pi\longrightarrow\Lambda
\longrightarrow\G\longrightarrow 1$, defined as follows
$$\Lambda=\left\{\;\lambda\in\mbox{Diff}(\ush)\;|\;\exists\;
\gamma\in\Gamma\;\mbox{s. t.}\;\pi(\lambda(u^{\#}))=\gamma\cdot
\pi(u^{\#})\;\forall\; u^{\#}\in\ush\;\right\}.$$ It is easy to verify 
that $\Lambda$
is a discrete group, that it acts on the manifold $\ush$ according 
to the assumptions
of Definition~\ref{model} and that $\ut/\G=\ush/\Lambda$.}\end{remark}
\begin{defn}[Submodel]{\rm
Consider a model $(\ut/\G,p,\ut)$ and let $W$ be an open subset of 
$\ut/\G$. We will
say that $W$ is a {\it submodel} of $(\ut/\G,p,\ut)$, 
if $(W,p,p^{-1}(W))$ defines a model by means of Remark~\ref{simplyc}.}
\end{defn}
\begin{defn}[Smooth mapping, diffeomorphisms of models]{\rm
Given two models $(\ut/\G, p, \ut)$ and $(\tilde{W}/\Delta,q,\tilde{W})$, 
a mapping
$f\,\colon\,\ut/\G\longrightarrow\tilde{W}/\Delta$ is said to be 
{\em smooth} if there
exists a smooth mapping $\tilde{f}\,\colon\,\ut\longrightarrow\tilde{W}$ 
such that $q\circ
\tilde{f}= f\circ p$; we will say that $\tilde{f}$ is a {\em lift} of $f$. 
We will say that the
smooth mapping $f$ is a {\em diffeomorphism of models} if it is 
bijective and if
the lift $\tilde{f}$ is a diffeomorphism.}
\end{defn}
If the mapping $\tilde{f}$ is a lift of a smooth mapping of models
$f\,\colon\,\ut/\G\longrightarrow \tilde{W}/\D$ so are the mappings
$\tilde{f}^{\gamma}(-)=\tilde{f}(\gamma\cdot -)$, for all elements 
$\gamma$ in $\G$ and
$^{\delta}\tilde{f}(-)=\delta\cdot\tilde{f}(-)$, for all elements 
$\delta$ in $\Delta$. If the
mapping $f$ is a diffeomorphism, then these are the only other possible 
lifts and the
groups $\G$ and $\D$ are isomorphic; for a proof see
\cite[orange and green lemmas]{p}.

Quasifolds are obtained by gluing together the models in the appropriate way:
\begin{defn}[Quasifold]
  {\rm A dimension $k$ {\em quasifold structure} on a topological space
    $M$ is the assignment of an {\em atlas}, or collection of {\em charts},
   ${\cal A}= \{\,(\ua,\fia,\uta/\ga)\,|\,\alpha\in A\,\}$ having the following properties:
\begin{enumerate}
\item[(1)] 
The collection $\{\,\ua\,|\,\alpha\in A\,\}$ is an open cover of $M$;
\item[(2)] For each index $\alpha$ in $A$ the space
$\uta/\ga$ defines a model, where the set $\uta$ is an open, connected, 
and simply
connected subset of the space $\R^k$, and the mapping $\fia$ is a homeomorphism of the
space $\uta/\ga$ onto the set $\ua$;
\item[(3)] For all indices $\alpha, \beta$ in $A$ such that
  $\ua\cap\ub\neq\emptyset$, the sets $\fia^{-1}(\ua\cap\ub)$ and
  $\fib^{-1}(\ua\cap\ub)$ are submodels of $\uta/\ga$ and $\utb/\gb$
  respectively and the mapping
  $$\gab=\fib^{-1}\circ\fia\,\colon\fia^{-1}(\ua\cap\ub)
  \longrightarrow\fib^{-1}(\ua\cap\ub)$$
  is a diffeomorphism of models. We will then say that the mapping
  $\gab$ is a {\em change of charts} and that the corresponding charts are
  {\em compatible};
\item[(4)] The atlas $\cal A$ is maximal, that is: if the triple
$(U,\phi,\ut/\G)$ satisfies property (2) and is compatible with all the 
charts in
$\cal A$, then $(U,\phi,\ut/\G)$ belongs to $\cal A$.
\end{enumerate}}
\end{defn}
\begin{remark}\label{gruppino}{\rm To each point $m\in M$ there corresponds a discrete group $\G_m$
defined as follows: take a chart $(\ua,\fia,\uta/\ga)$ around $m$, 
then $\G_m$ is the
isotropy group of $\ga$ at any point $\tilde{u}\in \ut_{\alpha}$ 
which projects down to $m$.
One can check that this definition does not depend on the choice of the
chart.}\end{remark}
\begin{defn}[Submodel in a quasifold]
\label{submodelofaquasifold}{\rm An open subset $W$ of $M$ is a 
{\it submodel in the quasifold M} if there exists a chart 
$(U,\phi,\tilde{U}/\G)$ of $M$ such that $W\subset U$ and 
$\phi^{-1}(W)$ is a submodel of 
$\tilde{U}/\G$.}\end{defn}
\begin{defn}[Smooth mapping, diffeomorphism]\label{diffeo}
{\rm Let $M$ and $N$ be two quasifolds. A continuous mapping
$f\,\colon\,M\longrightarrow N$ is said to be  {\em smooth}  if there
exists a chart $(\ua,\fia,\uta/\ga)$ around each point $m$ in the space 
$M$, a chart
$(\wa,\pia,\wta/\da)$ around the point $f(m)$, and a smooth mapping of models
$\fa\,\colon\,\vta/\ga\rightarrow\wta/\da$ such that $\pia\circ\fa=f\circ\fia$.
 If
$f$ is bijective, and if each $\fa$ is a diffeomorphism of models, 
we will say that
$f$ is a {\em diffeomorphism}.}
\end{defn}
\begin{defn}[Quasilattice, quasitorus]\label{quasilattice}{\rm Let $\d$ be a vector space of 
dimension $n$.
A {\em quasilattice}
in $\d$ is the $\Z$-span, $Q$, of a set of
$\R$-spanning vectors $X_1,\ldots,X_d\in\d$.
We call {\em quasitorus} of dimension $n$ the group and quasifold, covered 
by one chart, $D=\d/Q$.}
\end{defn}
Remark that $Q$ is a true lattice if and only if 
$\hbox{dim}_{\Q}\hbox{Span}_{\Q}\{X_1,\cdots,X_d\}=n$. In this case
$\d/Q$ is a torus.
\section{Stratifications by quasifolds}
We define the notion of space stratified by quasifolds in the
generality we need for our purposes. For the general definition of
stratification see \cite{GMcP1}. 
\begin{defn}\label{decomposition}{\rm
Let $M$ be a  topological space.
A {\it decomposition of $M$ by quasifolds} is a collection of
disjoint locally closed connected
quasifolds ${\cal T}_{F}$ ($F\in{\cal F}$),
called {\it pieces}, such that
\begin{enumerate}
\item[(1)] The set $\cal F$ is finite and partially ordered;
\item[(2)] $M=\bigcup_F{\cal T}_F$;
\item[(3)] ${\cal T}_F\cap{\overline{\cal T}}_{F'}\neq\emptyset$ if and only if
 ${\cal T}_F\subseteq{\overline {\cal T}}_{F'}$ if and only if $F\leq F'$.
\end{enumerate}}\end{defn}
We also require that ${\cal F}$ has a maximal element $F$ and that
the corresponding piece ${\cal T}_F$ is open and dense in $M$. We
call this piece the regular piece, the other pieces are called
singular. We will then say that $M$ is an {\em $n$-dimensional
compact space decomposed by quasifolds}, with $n$ the dimension of
the regular set.
\begin{remark}\label{cono}{\rm For the definition of 
stratification we need the following construction: let $L$ be
a compact space decomposed by quasifolds, we will call {\it cone over $L$},
denoted by $C(L)$, the space $[0,1)\times
L/\sim$, where two points $(t,l)$ and $(t',l')$ in $[0,1)\times L$
are equivalent if and only if $t=t'=0$. This space is itself a
space decomposed by quasifolds: for example when $L$ is a compact
quasifold the space $C(L)$ decomposes into two
pieces: one is the cone point, the other is given by the quasifold
$(0,1)\times L$. A further construction to be considered is the
following: let $t$ be a point in a quasifold $\cal T$,
$B\cong \tilde{B}/\Gamma$ 
a submodel in $\cal T$ containing $t$, 
and $L$ a compact space decomposed
by quasifolds. Notice first that
the decomposition of $L$ induces a decomposition of
the product $\tilde{B}\times\cono(L)$: to each piece $\cal L$ of $L$ there
corresponds the piece $\tilde{B}\times(0,1)\times{\cal L}$; to cover the
whole of $\tilde{B}\times\cono(L)$ we add a minimal piece, lying in the
closure of all other pieces, given by $\tilde{B}$ times the cone
point. Suppose, in addition, that $\G$ acts freely on $\tilde{B}$ and that
the space $L$ is
endowed with an
 action of $\Gamma$ that preserves the decomposition.
Then
the product $\tilde{B}\times\cono(L)$ is acted on by $\G$ and
the quotient $(\tilde{B}\times\cono(L))/\Gamma$ inherits the 
decomposition of $\tilde{B}\times\cono(L)$ in pieces. Moreover
the quotient $(\tilde{B}\times\cono(L))/\Gamma$ fibers over $B$ with fiber
$\cono(L)$. 
}\end{remark} 
A stratification is a decomposition that is
locally well behaved.
\begin{defn}\label{stratificazione}{\rm Let $M$ be an $n$-dimensional
compact  space decomposed by quasifolds, the decomposition of  $M$
is said to be a {\it stratification by quasifolds} if each
singular piece ${\cal T}$, called {\it stratum}, satisfies the
following conditions:
\begin{enumerate}
\item[(i)] let $r$ be the dimension of $\cal T$,
for every point $t\in{\cal T}$ there exist: 
an open neighborhood
$U$ of $t$ in $M$; 
a submodel $B\cong\tilde{B}/\Gamma$ in $\cal T$ containing 
$t$ and such that $\G$ acts freely on $\tilde{B}$;
an $(n-r-1)$-dimensional compact  space $L$ decomposed by
quasifolds, called the {\it link} of $t$;
an action of the group $\Gamma$ on $L$, 
preserving the decomposition of $L$ 
and such that the pieces of the induced decomposition of
 $\tilde{B}\times\cono(L)/\Gamma$ are quasifolds;
finally a homeomorphism
$h\,\colon\,(\tilde{B}\times \cono(L))/\Gamma\lorw U$ that
respects the decompositions and  
takes each piece of $(\tilde{B}\times\cono(L))/\Gamma$
diffeomorphically into the corresponding piece of $U$;
\item[(ii)] the decomposition of the link $L$ satisfies condition (i).
\end{enumerate}
}\end{defn}
The definition is recursive and, since the dimension of $L$ decreases
at each step, we end up, after a finite number of steps,
with links that are compact quasifolds.
\begin{remark}\label{quandobanale}{\rm Notice that, 
if 
the discrete groups $\G$'s are finite 
for any possible $F$, $t\in{\cal T}_F$ and $B$, then
the twisted products  $\tilde{B}\times\cono(L)/\Gamma$
become trivial and the singular strata turn out to be smooth manifolds, 
since $\G$'s act freeely. Therefore our stratification satisfies in this case 
the local triviality condition of the classical definition of 
stratification, morover, strata are smooth, with the only possile exception
of the principal stratum, that might be an orbifold.
We shall learn from the examples how the twisting discrete group
$\G$'s arise naturally from the contruction.
}\end{remark}
\section{The construction}
Let $\d$ be a real vector space of dimension $n$, and let $\Delta$
be a convex polytope of dimension $n$ in the dual space $\ddu$. We
want to associate to the polytope $\D$ a family of compact spaces
that are suitably stratified by symplectic quasifolds. We
construct these spaces as symplectic quotients, following the
procedure which was first introduced by Delzant in \cite{delzant}
and then extended to nonrational simple convex 
polytopes by Prato in \cite{p}.
Write the polytope as
\begin{equation}\label{polydecomp}
\D=\bigcap_{j=1}^d\{\;\mu\in\ddu\;|\;\langle\mu,X_j\rangle\geq\lambda_j\;\}
\end{equation}
for some elements $\xd$ in the vector space $\d$  
and some real
numbers $\ld$. 
Let $Q$ be a quasilattice in the space $\d$
containing the elements $X_j$ (for example the one that is
generated by these elements, namely $\hbox{Span}_{\Z}\{X_1,\dots,X_d\}$) 
and let $\{\ed\}$ denote the standard
basis of $\rd$; consider the surjective linear mapping
\begin{equation}\label{projection}
\begin{array}{cccc}\label{pi}
\pi \,\colon\,& \R^d & \lorw & \d\\
    &   e_j& \longmapsto & X_j.
\end{array}
\end{equation}
Consider the $n$-dimensional quasitorus $\d/Q$. The mapping $\pi$
induces a group homomorphism,
\begin{equation}\label{defdienne}
\Pi
\,\colon\, T^d=\rd/\Z^d\lorw \d/Q.
\end{equation}
We define $N$ to be the kernel of the mapping $\Pi$.
The mapping $\Pi$ defines an isomorphism
\begin{equation}\label{qtiso}
T^d/N\longrightarrow \d/Q.
\end{equation}
We construct a moment mapping for the Hamiltonian action of $N$ on
$\C^d$.
 Consider the mapping $$\Upsilon(\vz)=\sum_{j=1}^d
(|z_j|^2+\lambda_j)e_j^*,$$ where the $\lambda_j$'s are given in
(\ref{polydecomp}) and are uniquely determined by our choice of
inward pointing normals to codimension $1$ faces. The mapping $\Upsilon$ is a moment mapping for the
standard action of $T^d$ on $\C^d$. Consider now the subgroup
$N\subset T^d$ and the corresponding inclusion of Lie algebras
$\iota\,\colon\,\n\rightarrow\R^d$. The mapping
$\Psi\,\colon\,\C^d\rightarrow \n^*$ given by $\Psi={\iota}^*\circ
\Upsilon$ is a moment mapping for the induced action of $N$ on $\cd$. We
want to prove that the quotient $M=\zset/N$, endowed with the
quotient topology, is a space stratified by  quasifolds. Notice
that, by (\ref{defdienne}), the group $N$ is not 
closed in $T^d$ unless $Q$ is an honest lattice, 
moreover to each $\D$ there corresponds a whole
family of quotients, given by all possible choices of normal
vectors and of quasilattices $Q$ containing these vectors.
Let us remind that
the polytope $\D$ is said to be {\em rational} if there exist a lattice
$L$ and a a choice of normals $X_j$ such that $X_j\in L$ for $j=1,\dots,d$.
Nonrational polytopes are for example the regular pentagon and Penrose's kite,
whilst triangles are all rational. However  
we can associate to a rational polytope a ``nonrational'' space by making a 
nonrational choice of data attached to it. More precisely
a choice of normals and quasilattice $Q$ is said to be {\it rational} 
if  
$Q$ is a true 
lattice, it is {\it nonrational} if the chosen
quasilattice $Q$ is not a lattice.

In our setting, in which the polytope can be nonsimple,
the zero set $\zset$ is not in general a smooth submanifold of
$\R^{2d}$. Nonsimplicity of the polytope is responsible, as in
the rational case, for the decomposition in strata of the
quotient.

To define the decomposition of $M$ in pieces we start by giving
some further definitions on the polytope $\D$.
Let us consider the open faces of $\Delta$. They can be described
as follows. For each such face $F$ there exists a 
subset $I_F\subset\{1,\ldots,d\}$ such that
\begin{equation}\label{facce}
F=\{\,\mu\in\Delta\;|\;\langle\mu,X_j\rangle=\lambda_j\;
\hbox{ if and only if}\; j\in I_F\,\}.
\end{equation}
The $n$-dimensional open face of $\D$, which we denote
by $\hbox{Int}(\D)$, corresponds to the empty subset.
A partial order on the set of all faces of $\D$ is defined by
setting $F\leq F'$ (we say $F$ contained in $F'$) if
$F\subseteq\overline{F'}$. The polytope $\Delta$ is the disjoint
union of its faces. Let $r_F=\hbox{card}(I_F)$; we have the
following definitions:
\begin{defn}{\rm A $p$-dimensional face $F$ of the polytope is
said to be {\em singular} if $r_F>n-p$,
{\em nonsingular} if $r_F=n-p$.}
\end{defn}
\begin{remark}{\rm Let $F$ be a $p$-dimensional singular face in $\d^*$, 
then $p<n-2$.
For example: a polytope in $(\R^2)^*$ is simple; the singular faces of a 
nonsimple polytope in $(\R^3)^*$ must be $0$-dimensional.
}\end{remark}
The following Proposition is an adaptation of an analogous results
in \cite{g,p}, to which we refer the reader for further details.
\begin{prop}\label{phi} The $n$-dimensional quasitorus $D=\d/Q$ acts
continuously and effectively on the topological space $M=\zset/N$. Moreover $M$ is
compact and a continuous mapping $\Phi\,\colon\,M\lorw\d^*$ is
defined such that $\Phi(M)=\D$.\end{prop} \proof
Consider the
exact sequence
\begin{equation}
\label{exactsequence}
0\lorw\d^*\stackrel{\pi^*}{\lorw}(\R^d)^*\stackrel{{\iota}^*}{\lorw}(\n)^*\lorw0.
\end{equation}
By (\ref{exactsequence}) we have that
the mapping $(\pi^*)^{-1}\circ \Upsilon$ induces a continuous mapping from the 
quotient $M$ to $\d^*$,
we call this mapping $\Phi$. Now notice, by making use of the explicit expression of $\Upsilon$, that $\vz\in\zset$ if and only if
\begin{equation}
\label{formula} 
|z_j|^2=\langle \Phi([\vz]),X_j\rangle-\lambda_j,\quad\quad j=1,\ldots, d.
\end{equation}
This implies that $\Phi(M)=\D$. Moreover properness of $\Upsilon$ implies
that $M$ is compact. A continuous action 
\begin{equation}\label{action}
\tau\colon D\times M\lorw M
\end{equation}
is defined via
the isomorphism (\ref{qtiso}).
This action is free on $\Phi^{-1}(\hbox{Int}(\D))$.
\qed
\begin{remark}{\rm 
An immediate consequence of the arguments used in the proof of 
Proposition~\ref{phi} is that the mapping $\Phi$ induces a homeomorphism
from the topological quotient $M/D$ onto $\D$.
}\end{remark}
We introduce now two examples of nonsimple convex polytopes that have served 
as models for the whole construction. They are both rational polytopes
with respect to the integer lattice $\Z^n\subset\R^n$, that is they 
both admit choices of normals contained in the integer lattice.
Hence they admit rational and nonrational choices; 
the corresponding families of spaces
include therefore spaces stratified by smooth manifolds/orbifolds 
and spaces stratified by quasifolds  (cf. Remark~\ref{choices}). 
Our first example is a polytope in $(\R^3)^*$, the second one is a polytope in 
$(\R^4)^*$. This last 
example allows us to illustrate the features of our spaces to a greater extent,
since it has singular faces of positive dimension. 
Both examples will be resumed in the last section.
\begin{example}[Figure 1]\label{piramide}{\rm In $(\R^3)^*$ 
consider the pyramid given by the convex hull of 
$\nu=(0,0,1)$ (the apex), $\mu_1=(1,0,0),\mu_2=(1,1,0),\mu_3=(0,1,0),
\mu_4=(0,0,0)$. It is a nonsimple polytope with only a singular face: 
the apex $\nu$. 
We make the following choice of inward pointing normals
\begin{equation}
\begin{array}{lll}
\langle \mu_1,\mu_2,\nu\rangle&X_1=(-1,0,-1)&\lambda_1=-1\\
\langle \mu_2,\mu_3,\nu\rangle&X_2=(0,-p_2,-p_2)&\lambda_2=-p_2\\
\langle \mu_3,\mu_4,\nu\rangle&X_3=(1,0,0)&\lambda_3=0\\
\langle \mu_4,\mu_1,\nu\rangle&X_4=(0,1,0)&\lambda_4=0\\
\langle \mu_1,\mu_2,\mu_3,\mu_4\rangle&X_5=(0,0,p_5)&\lambda_5=0\end{array}
\end{equation}
with $p_j\in\R_{>0}$ for $j=2,5$.
The $j$-th item of the first column lists the vertices 
contained in the $1$-codimensional affine space $<\mu,X_j>=\lambda_j$.
We choose the quasilattice $Q$ to be the one generated by $X_1,\dots,X_5$. 
}\end{example}
\begin{figure}[h]
\begin{center}
\input piramide.pstex_t
\end{center}
\caption{Ex.~\ref{piramide}}
\end{figure}
\begin{example}[Figure 2]\label{tenda}{\rm \hspace{7pt} In $(\R^4)^*$ consider the convex hull
of \linebreak $\nu_1=(1,-1,0,0)$, $\nu_2=(0,0,1,-1)$  and 
$\mu_1=(1,0,0,0),$ $\mu_2=(0,1,0,0),$ $\mu_3=(0,0,1,0),$ $\mu_4=(0,0,0,1)$. 
The resulting polytope is nonsimple, 
with nine $3$-dimensional faces.
It can be thought of as the 4-simplex in which the origin has been substituted
by the edge $\nu_1\nu_2$.
The singular faces are all of its 
vertices and some of its edges.
We make the following choice of inward pointing normals
\begin{equation}
\begin{array}{lll}
\langle \nu_1,\nu_2,\mu_1,\mu_2\rangle&X_1=(0,0,p_1,p_1)&\lambda_1=0\\
\langle \nu_1,\nu_2,\mu_3,\mu_4\rangle&X_2=(1,1,0,0)&\lambda_2=0\\
\langle \nu_1,\nu_2,\mu_1,\mu_3\rangle&X_3=(-1,0,-1,0)&\lambda_3=-1\\
\langle \nu_1,\nu_2,\mu_2,\mu_4\rangle&X_4=(2,1,2,1)&\lambda_4=1\\
\langle \nu_2,\mu_1,\mu_2,\mu_3\rangle&X_5=(-p_5,-p_5,-p_5,0)&\lambda_5=-p_5\\
\langle \nu_1,\mu_1,\mu_2,\mu_4\rangle&X_6=(0,0,1,0)&\lambda_6=0\\
\langle \nu_1,\mu_1,\mu_3,\mu_4\rangle&X_7=(-1,0,-1,-1)&\lambda_7=-1\\
\langle \nu_2,\mu_2,\mu_3,\mu_4\rangle&X_8=(p_8,0,0,0)&\lambda_8=0\\
\langle \mu_1,\mu_2,\mu_3,\mu_4\rangle&X_9=(-1,-1,-1,-1)&\lambda_9=-1
\end{array}
\end{equation}
with $p_j\in\R_{>0}$ for $j=1,5,8$.
The $j$-th item of the first column lists the vertices 
contained in the $1$-codimensional affine space $<\mu,X_j>=\lambda_j$.
Thus we find that $r_{\eta}=6$ for each vertex $\eta$ of the polytope.
The singular edges are those with
$r_F=4$ (marked in red in the picture).
For example: $I_{\nu_1}=\{1,2,3,4,6,7\}, r_{\nu_1}=6$,
 $I_{\nu_2}=\{1,2,3,4,5,8\}, r_{\nu_2}=6$,
 $I_{\nu_1\nu_2}=\{1,2,3,4,\}, r_{\nu_1\nu_2}=4$,
where $\nu_1\nu_2$ denotes the edge joining $\nu_1$ to $\nu_2$.
We choose the quasilattice $Q$ to be the one generated by $X_1,\dots,X_9$. 
}
\end{example}
\begin{figure}[h]
\begin{center}
\input tenda.pstex_t
\end{center}
\caption{Ex.~\ref{tenda}}
\end{figure}
\section{Notation}
We gather in this section the necessary notation  to proceed
with the proofs of the main results.
Let $\K$ be one of the following sets $\C,\C^*,\R,\R^*,\R_{\geq0},\R_{>0}$,
all of them considered naturally immersed in $\C$.
Let $J$ be a subset of $\{1,\dots,d\}$ and let $J^c$ be its complement.
We denote by
$$\K^J=\{\zd\in\C^d\;|\;z_j\in\K\;\;\hbox{if}\;\; j\in J,\;
z_j=0\;\;\hbox{if}\;\; j\notin J
\}.$$
For example if $J=\emptyset$ then 
$(\K)^{J^c}=\{\zd\in\C^d\;|\;z_j\in\K\}=\K^d$.
We have $\K^d=\K^J\times(\K)^{J^c}$.
Let $\vz\in\K^d$, we denote by $\vz_J$ its projection onto the factor
$\K^J$.
The stabilizer of the $T^d$-action at any point in ${(\C^*)}^{J^c}$ 
is the torus
$$T^{J}=\{\,(a_1,\dots,a_d)\in T^d\;|\;a_j=1\;\;\hbox{if}\;\;j\notin J\,\}$$
of Lie algebra $\R^J$.
Let $F$ be a $p$-dimensional face and let $I_F$ be the corresponding
set of indices. To lighten the notation we omit the $I$ and write 
$T^F$ for $T^{I_F}$, $\R^F$ for $\R^{I_F}$, ${(\C^*)}^{F^c}$ for ${(\C^*)}^{I_F^c}$, etc. 
Moreover we set
\begin{equation}\label{enneeffedef}
N^F=N\cap T^F\end{equation} and $$\n^F=\n\cap \R^{F}$$ its Lie algebra.
Let $\d_F=\pi(\R^{F})$, where $\pi:\R^d\lorw\d$ is the projection defined 
in (\ref{projection}); notice that $\d_F\cap Q$ is a quasilattice.
The subgroup $N^F$ of $N$ has dimension $(r_F-n+p)$, moreover
$T^ F/N^F\cong\d_F/(\d_F\cap Q)$.
Now let $\mu$ be a vertex of $\D$.
Define ${\cal I}_{\mu}$ to be the set of subsets $I$ of $I_{\mu}$ such that
the set $\{X_j\;|\;j\in I\}$ is 
a basis of  $\d$.
If $\mu$ is non singular then ${\cal I}_{\mu}$ contains just the element
$I_{\mu}$.
We denote by 
$${\cal I}=\bigcup_{(\mu\;\hbox{vertex of}\;\D)}{\cal I}_{\mu}.$$
For each $I\in{\cal I}$ we have
that the group $N\cap T^I$ is discrete, since $\pi(\R^I)$ is
$n$-dimensional. 
We set
\begin{equation}\label{gammai}
\Gamma_I=N\cap T^I,\quad\hbox{for}\quad I\in{\cal I}.
\end{equation}
Now let $F$ be a singular face of $\D$ of dimension $p>0$ and
let $\mu$ be a vertex contained in $F$, hence $I_F\subset I_{\mu}$.
Take $I\in{\cal I}_{\mu}$, such that $\hbox{card}(I\cap I_F)=(n-p)$.
We define the discrete group
\begin{equation}\label{gammaeffei}
\G_{I\cap I_F}=N^F\cap T^{I\cap I_F} .
\end{equation}
Consider now the diagrams
\begin{equation}\label{primamappa}
\Gamma_I\hookrightarrow T^{I\cap I_F}\times T^{I\setminus (I\cap I_F)}
\longrightarrow T^{I\cap I_F}
\end{equation}
\begin{equation}\label{secondamappa}
\Gamma_I\hookrightarrow T^{I\cap I_F}\times T^{I\setminus (I\cap I_F)}
\longrightarrow T^{I\setminus(I\cap I_F)}.
\end{equation}
Define $\check{\Gamma}^I_{I\cap I_F}$ and 
$\check{\Gamma}_{I\setminus(I\cap I_F)}$
to be the images of the composition of
mappings (\ref{primamappa}) and (\ref{secondamappa}) respectively. 
The groups $\check{\Gamma}^I_{I\cap I_F}$ and 
$\check{\Gamma}_{I\setminus(I\cap I_F)}$
are discrete groups such that
$\Gamma_I=\check{\Gamma}^I_{I\cap I_F}\times 
\check{\Gamma}_{I\setminus(I\cap I_F)}$. Moreover
the following group exact sequence is defined: 
\begin{equation}
1\hookrightarrow
\G_{I\setminus I\cap I_F}\lorw\check{\Gamma}_{I\setminus(I\cap I_F)}
\stackrel{\ell_{I_F,I}}
{\lorw}\;\check{\Gamma}^I_{I\cap I_F}/\Gamma_{I\cap I_F}\lorw1.
\end{equation}
\begin{remark}\label{ifhonestlattice}{\rm
The discrete groups just defined have a crucial role in the sequel.
Notice that they need not to be finite, but they are so 
whenever $Q$ is an honest lattice. 
They are trivial if, for example, we can choose inward pointing normals
$\{X_1,\dots,X_d\}$ such that:
i) the quasilattice generated by $\{X_1,\dots,X_d\}$ is a lattice $L$
and we choose $Q=L$;
ii) for each
$I\in{\cal I}$ the set $\{X_j,\quad j\in I\}$ is a basis for $Q$.
In particular a Delzant polytope in a lattice $L$ realizes the above 
condition if 
the $X_j$'s are taken to be primitive in $L$.
}
\end{remark}
\section{The stratification}
We are now ready to define the decomposition of $M$, the 
singular pieces
are given by:
$${\cal T}_F=\Phi^{-1}(F)\quad\hbox{with}\;F\;\hbox{singular face};$$
the regular piece, which contains $\Phi^{-1}(\hbox{Int}(\Delta))$,
 is given by 
$${\cal T}_{\D}=\cup_{F\;\hbox{nonsing.}}
\Phi^{-1}(F).$$ 
\begin{remark}\label{stratidizset}{\rm It is important to point out 
that for any face $F$ of $\D$ the
set $\Phi^{-1}(F)$ is non-empty and is precisely given by
$(\zset\cap(\C^*)^{{F}^c})/N$; this follows easily
from the proof of
Proposition~\ref{phi} 
(for further details cf. \cite{g}).}\end{remark}
Remark~\ref{stratidizset} allows us
to characterize the pieces ${\cal T}_F$ in the standard way, by
the isotropy group attached to each of them.
\begin{remark}\label{stratidizset2}{\rm \hspace{2mm}
The regular piece ${\cal T}_{\D}$ is given by the quotient 
$\zset_{\D}/N$, where \linebreak 
$\zset_{\D}=\cup_{F\;\hbox{nonsing.}}(\zset\cap (\C^*)^{F^c}).$
Therefore the stabilizer of $N$ at each point of $\zset_{\D}$ is discrete
and $\zset_{\D}$ is precisely
 the set of regular points of 
the moment mapping
$\Psi$, while 
 the stabilizer 
of $N$ corresponding to each singular piece ${\cal T}_F$ is precisely   
the $(r_F-n+p)$-dimensional group
$N^F$ defined
in (\ref{enneeffedef}).
}
\end{remark}
\begin{thm}[Quasifold structure of strata]\label{stratiquasifold}
The subset ${\cal T}_F$ of $M$ corresponding to each $p$-dimensional
singular face of $\D$ is a $2p$-dimensional quasifold.
The subset ${\cal T}_{\D}$
is a $2n$-dimensional quasifold. These subsets give a
decomposition by quasifolds of $M$.
\end{thm}
Before giving the proof of this theorem we prove  
two key Lemmas, they ensure 
that the behavior
of the group
$N$ and its subgroups $N^F$ does not differ much
from that of a torus; noncompactness can always be 
concentrated within a discrete group. 
The first Lemma is the real version of \cite[Lemma 2.3]{cx}, we
recall the proof for completeness.
\begin{lemma}[The group $N$]  
\label{torustrick} Let $\mu$ be a vertex of the
polytope $\D$ and let $I\in{\cal I}_{\mu}$.
Then we have
that\\ {\rm (i)} $T^d/T^I\cong N/\Gamma_{I}$;\\ {\rm
(ii)} $N=\Gamma_{I}\exp{(\n)}$;\\ {\rm (iii)} given any
complement $\b$ of $\R^I$ in $\R^d$, we have that
$$\n=\{\,Y-\pi_{I}^{-1}(\pi(Y))\;|\;Y\in\b\,\}.$$
\end{lemma}
\proof (i) Consider the group homomorphism 
$$
\begin{array}{cccc}
\lambda_{I} \,\colon\,& N&\lorw&T^d/T^I\\
&n&\longmapsto&[n].
\end{array}
$$ Since $\n$ and $\R^I$ are complementary
 $\lambda_I$ is surjective. The kernel of
$\lambda_I$ is given by $\G_I$, therefore $\lambda_I$ induces an isomorphism
$T^d/T^I\cong N/\Gamma_I$.\\(ii) Every element in $N$ can be written
in the form $\exp{(X)}$, where $X\in\R^d$ is such that 
$\pi(X)\in Q$. Write now
$X=X-\pi_{I}^{-1}(\pi(X))+\pi_{I}^{-1}(\pi(X))$; it is easy to check that
$X-\pi_{I}^{-1}(\pi(X))\in\n$, and that
$\exp{(\pi_{I}^{-1}(\pi(X)))}\in\G_{I}$. The group $\G_{I}\cap\exp{(\n)}$ is
not necessarily trivial, so the decomposition is not necessarily unique.
\\ (iii) Every element of the form $Y-\pi_{I}^{-1}(\pi(Y))$, with
$Y\in{\b}$, clearly belongs to $\n$. 
Conversely, write every element $V\in\n$ as
$V=X+Y$ according to the decomposition $\R^d=\R^I\oplus\b$, and notice that
$\pi(V)=0$ implies that $X=-\pi_{I}^{-1}(\pi(Y))$. \qed
\begin{lemma}[The group $N^F$]\label{torusspliteffe} 
Let $F$ be a singular face of $\D$ and let
 $\mu$ be a vertex contained in $F$. For any choice
of $I\in{\cal I}_{\mu}$ such that
$\hbox{card}(I\cap I_F)=(n-p)$ we have:\\
{\rm (i)} $T^F/T^{(I\cap I_F)}\cong N^F/\G_{I\cap I_F}$;\\
{\rm (ii)} $N^F=\G_{I\cap I_F}\exp(\n^F)$;\\
{\rm (iii)} $\n^F=\{Y-\pi_{I\cap I_F}^{-1}(\pi_F(Y))
\;|\;Y\in\R^{I_F\setminus I\cap I_F}\}$. 
\end{lemma}
\proof See the proof of Lemma~\ref{torustrick}. Notice that also in this
case the intersection $\G_{I\cap I_F}\cap\exp(\n^F)$ may not be trivial, 
so the 
decomposition need not to be unique.
\qed
\begin{remark}\label{basefiltrazioneenne}{\rm 
We exhibit here, by means of Lemmas~\ref{torustrick} and \ref{torusspliteffe},
an explicit basis of $\n$, which is very
useful for explicit computations. Consider the 
flag in $\n$ 
\begin{equation}\label{filtrazione}
\n^F\subset\n^{\mu}\subset\n.
\end{equation}
Let $A_I=(a^I_{ij})_{
i\in I
}\in M_{n\times d}$ the matrix of the projection
$\pi:\R^d\lorw\d$ with respect to the standard basis
of $\R^d$ and the basis
$\{X_h,\;h\in I\}$ of $\d$. When clear from the context we will omit the 
$I$'s and 
write simply 
$A=(a_{ij})$.
A basis of $\n$ adapted to the flag (\ref{filtrazione}) is given by the 
following vectors:
\begin{equation}\label{baseenneeffe}
e_k-\sum_{h\in I\cap I_F}a_{hk}e_{h},\;\quad k\in I_F\setminus I\cap I_F;
\end{equation}
\begin{equation}\label{baseenne}
e_l-\sum_{h\in I}a_{hl}e_{h},\;\quad l\in I_{\mu}\setminus (I\cup I_F);
\quad\quad\quad
e_r-\sum_{h\in I}a_{hr}e_{j_h},\;\quad r\notin I_{\mu}. 
\end{equation}}
\end{remark}
{\bf Proof of Theorem~\ref{stratiquasifold}} We have already observed in 
Remarks~\ref{stratidizset} and \ref{stratidizset2} that the regular 
and singular pieces
are quotients of suitable subsets of $\zset$ by the group $N$.
In order to 
construct, for each piece, a quasifold atlas,  
we construct local slices for the corresponding subset of $\zset$.
This leads to the construction of local models.
The crucial point
is that the natural bijective mapping, from each local model thus obtained 
into the piece in consideration, is closed (cf. Step (I,d)).
We prove the statement separately 
for the regular and singular pieces and divide the proof in steps,
each with a title, in order to simplify the
exposition and make it easier to refer to parts of the proof.

{\it Part I: The regular stratum.} 
The fact that ${\cal T}_{\D}$ is a $2n$-dimensional symplectic quasifold
acted on by $D$ descends from the general result \cite[Thm~3.1]{p}. 
We need to give here an explicit proof.

{\it (I,a): Construction of local models}:
In order to prove that ${\cal T}_{\D}$
is a $2n$-dimensional quasifold we construct a collection of charts
covering ${\cal T}_{\D}$.
For each $I\in{\cal I}$ consider
the open subset $\uci$ of $\regpsi$ defined by 
$\uci=\zset\cap(\C^I\times(\C^*)^{I^c})$.
The group $N$ acts on each $\uci$ with discrete stabilizer and the open
subsets $\uci$ cover the whole regular set $\zset_{\D}$. 
Therefore
the quotients $\uci/N$, that we denote by $U_I$, 
give an open covering of 
${\cal T}_{\D}$.
Consider now a  vertex $\mu$ and an $I\in{\cal I}_{\mu}$. Let
$\G_I$ be the discrete group defined in (\ref{gammai}). 
We want to prove that
there exist an open subset $\uti\subset\C^n$
and a mapping $\phi_I:\uti/\Gamma_I\lorw U_I$ 
such that
$\uti/\Gamma_I$ is a model, in the sense on Remark~\ref{simplyc}, and 
the mapping $\phi_I$ is a homeomorphism.
The components of the moment mapping $\Psi$ 
with respect to the basis dual to the adapted basis given in 
(\ref{baseenneeffe},\ref{baseenne}) are as follows:
the first $r_{\mu}-n$ components are given by
\begin{equation}\label{primecomponenti}
-\sum_{h\in I}a_{hk}(|z_{h}|^2+\lambda_{h})+(|z_k|^2+\lambda_k)
\end{equation}
with $k\in I_{\mu}\setminus I$.
Since for each $k\in I_{\mu}\setminus I$ we have  
\begin{equation}\label{iperpiani}
\lambda_k=<\mu,X_k>=<\mu,\sum_{h\in I}a_{hk}X_{h}>=
\sum_{h\in I} a_{hk}<\mu,X_{h}>=\sum_{h\in I} a_{hk}\lambda_{h}
\end{equation}
the expression (\ref{primecomponenti}) reduces to
\begin{equation}\label{psimu}
-\sum_{h\in I}a_{hk}|z_{h}|^2+ |z_k|^2
\end{equation}
with $k\in I_{\mu}\setminus I$.
The remaining $d-r_{\mu}$ components are given by
\begin{equation}\label{psi}
-\sum_{h\in I}a_{hr}(|z_{h}|^2+\lambda_{h})+(|z_r|^2+\lambda_r)
\end{equation}
where $r\notin I_{\mu}$ and, by (\ref{iperpiani}), 
$(\sum_{h\in I}a_{hr}\lambda_{h})-\lambda_r>0$.
We define the set $\uti$ by
$$\begin{array}{rcll}\uti&=&
\{\vu\in\C^I\;|&\;\sum_{h\in I}a_{hk}|u_{h}|^2>0\;\;\hbox{for}\;\;
k\in I_{\mu}\setminus I,\quad\quad\quad\quad\\
&&&\sum_{h\in I}a_{hr}|u_{h}|^2>-((\sum_{h\in I}a_{hr}\lambda_{h})-\lambda_r)\;\;\hbox{for}\;\;r\notin I_{\mu}\}.\end{array}$$
The set $\uti$ is a nonempty open subset of $\C^I\cong\C^n$, star-shaped
with respect to the origin,
the quotient $\uti/\G_I$ is a model, in the sense on Remark~\ref{simplyc}.
We construct now the homeomorphism 
$\phi_I:\uti/\Gamma_I\lorw U_I$.
Consider the mapping 
\begin{equation}\label{mappaeffei}
\begin{array}{cccc}
F_I\,\colon& \uti&\lorw &(\R_{>0})^{I^c}\\
\end{array}
\end{equation}
defined as follows 
$$(F_I)_k(\vu)=0\quad k\in I$$
$$(F_I)_l(\vu)=\sqrt{\sum_{h\in I}a_{hl}|u_{h}|^2}\quad
l\in I_{\mu}\setminus I$$
and 
$$(F_I)_r(\vu)=\sqrt{\sum_{h\in I}a_{hr}|u_{h}|^2+
(\sum_{h\in I}a_{hr}\lambda_{h}-\lambda_r)}\quad r\notin I_{\mu}.$$
Now 
define the following mapping
$$
\begin{array}{ccc}
 \uti&\lorw&U_I\\
\vu&\longmapsto&[\vu+F_I(\vu)].
\end{array}
$$ 
We want to prove that the induced mapping $\phi_I$ from $\uti/\Gamma_I$ 
to $U_I$
is a homeomorphism.
Observe first that the induced mapping $\phi_I$ is  continuous, and
injective by 
definition of $\Gamma_I$.
It is also surjective, this can be proved as follows:

{\it (I,b) Surjective mapping}:
In order to prove surjectivity of $\phi_I$ take any element
$\vu+\vw\in\uci$ with $\vu\in\C^I,\;\vw\in(\C^*)^{I^c}$.
By Lemma~\ref{torustrick} we can choose an element $a\in N$ such that
$a\cdot(\vu+\vw)=\vu'+\vw'$ with $\vu'\in\C^I$ and 
$\vw'\in(\R_{>0})^{I^c}$.
Therefore $\vw'=F_I(\vu')$ and $\phi_I$ is surjective.
In other words 
the image of the mapping 
$$\begin{array}{cccc}
\tilde{\phi}_I\colon&\uti&\lorw&\uci\\
&\vu&\longmapsto&\vu+F_I(\vu),
\end{array}
$$
is a slice whose saturation is exactly $\uci$.

{\it (I,c) Closed mapping}: In order to prove that $\phi_I$ is closed, 
we need to check that $N(\tilde{\phi}_I(C))$
is closed for each $\Gamma_I$-invariant closed subset $C$ of $\uti$.
Let $a_m(\vu'_m+F_I(\vu'_m))$ 
be a sequence in $N(\tilde{\phi}_I(C))$ converging to
$\vx$, which necessarily  belongs 
to $\zset$. We prove that $\vx$ lies
in $N(\tilde{\phi}_I(C))$.
By
Lemma~\ref{torustrick}--the key fact here--the 
elements $a_m$ of $N$ can be decomposed as $g_m\exp(-\pi_I^{-1}(\pi(Y'_m)))
\exp(Y'_m)$, for suitable elements $Y'_m\in\R^{I^c}$ and $g_m\in\G_I$.
Remark now that the sequence $\exp(Y'_m)$ is in the torus $T^{I^c}$, 
therefore it does admit a subsequence, that we call in the same
way, converging to an element of the form $\exp(Y)$, with $Y\in\R^{I^c}$.
We can now choose a sequence $Y_m\in\R^{I^c}$ such that 
$\exp{Y_m}=\exp{Y'_m}$, for each $m$, with $Y_m$ converging to $Y$.
This leads to $a_m=k_mb_m$,  with $k_m\in\G_I$ and 
$b_m=\exp(-\pi_I^{-1}(\pi(Y_m)))
\exp(Y_m)$ converging to $b=\exp(-\pi_I^{-1}(\pi(Y)))
\exp(Y)$ by continuity. Therefore 
$a_m(\vu'_m+F_I(\vu'_m))$ becomes
$b_m(\vu_m+F_I(\vu_m))$, where $\vu_m=k_m\vu'_m$--here we use the
invariance of $F_I$ under the action of $\G_I$.
By continuity of the $T^d$-action this implies that 
$(\vu_m+F_I(\vu_m))\lorw b^{-1}\vx$, in particular the sequence
$\vu_m$ is convergent and its limit $\vw$ is in $C$, since $C$ is closed.
Therefore we can deduce that $\vx=b(\vw+F_I(\vw))$.
 
{\it (I,d) Universal covering}: Now we compute the fundamental group of $\uti$.
Denote by 
$$\begin{array}{r}
C_I=\{\vrho\in
(\R_{\geq0})^I\;|\;
\sum_{h\in I}a_{hk}\rho_{h}>0\;\hbox{for}\;
k\in I_{\mu}\setminus I,\\
\sum_{h\in I}a_{hr}\rho_{h}>-(\sum_{h\in I}a_{hr}\lambda_{h}-\lambda_r)\;\hbox{for}\;r\notin I_{\mu}\}.\end{array}$$
The set $C_I$ is an intersection of half-spaces, it is therefore convex.
Denote by $\{\rho_{h}=0\}$ the coordinate hyperplane 
$\{\vrho\in
\R^I\;|\;\rho_{h}=0\}$. Let 
$I_*=\{h\in I\;|\;\{\rho_{h}=0\}\cap C_I=\emptyset\}$ and let
$\ell=\hbox{card}(I_*)$, then $\pi_1(\uti)=\Z^{\ell}$. When $\ell=0$, the 
quotient
$\uti/\G_I$ is a model, when
$\ell>0$ we construct a chart by taking the universal
covering $U^{\sharp}_I$ of $\uti$ and the discrete group $\Lambda_I$, 
extension of $\Gamma_I$ by $\pi_1(\uti)$, as explained in
Remark~\ref{simplyc}.
We thus obtain a model homeomorphic to $U_I$.

{\it (I,e) Change of charts}: 
To prove that the charts constructed above 
are compatible
we need to check that the changes of charts are 
diffeomorphisms of models, more precisely:
consider two subsets $I$ and $J$ in ${\cal I}$.
Suppose that the corresponding
charts $U_I$ and $U_J$ have nonempty intersection.
We want to prove that the mapping
$$\phi_J^{-1}\circ\phi_I\,:\phi_I^{-1}(U_I\cap U_J)\lorw
\phi_J^{-1}(U_I\cap U_J)$$
is a diffeomorphism of models.   
For simplicity we consider the case in which $\uti$ and $\tilde{U}_J$
are both simply connected. We adapt the proof given in \cite[Thm~2.2]{cx}.
Let $$\tilde{W_I}=
\uti\cap\left(\C^{I\cap J}\times{\C^*}^{I\setminus(I\cap J)}
\right)$$
and
$$\tilde{W}_J=
\utj\cap\left(\C^{I\cap J}\times{\C^*}^{J
\setminus(I\cap J)}\right).$$
Then $W_I=\tilde{W}_I/\G_I$ is exactly $\phi_I^{-1}(U_I\cap U_J)$
and $W_J=\tilde{W}_J/\G_J$ is exactly $\phi_J^{-1}(U_I\cap U_J)$, so they
are submodels of $\uti/\G_I$ and $\utj/\G_J$ as required
by the definition.
In order to have simply connected open sets we pass to the universal
coverings $W^{\sharp}_I$ and $W^{\sharp}_J$ of $\tilde{W}_I$ and $\tilde{W}_J$
respectively.
We have
$$W^{\sharp}_I=\{(\vu,\vrho,\vtheta)\in\C^{I\cap J}\times(\R_{>0})^{I\setminus(I\cap J)}\times
\R^{I\setminus(I\cap J)}\;|\;(\vu,\sqrt{\vrho}\exp(\vtheta))\in\uti\}
$$
and 
$$W^{\sharp}_J=\{(\vu,\vrho,\vtheta)\in
\C^{I\cap J}\times(\R_{>0})^{J\setminus(I\cap J)}\times
\R^{J\setminus(I\cap J)}\;|\;(\vu,\sqrt{\vrho}\exp(\vtheta))\in\utj\}$$
where $(\sqrt{\vrho}\exp\vtheta)_{\{j\}}=\sqrt{\rho_j}\exp(\vtheta_j)$.
These are simply connected open sets acted on by the discrete groups
$$\Lambda_I=\{
(\exp X,Y)\;|\;X\in\R^{I\cap J}, Y\in\R^{I\setminus(I\cap J)},\;\pi(X+Y)\in Q\}$$
and
$$\Lambda_J=\{
(\exp X,Y)\;|\;X\in\R^{I\cap J}, Y\in\R^{J\setminus(I\cap J)},\;\pi(X+Y)\in Q\}$$
in the following manner:
$$
\begin{array}{lcrcc}
(\Lambda_I&,&W^{\sharp}_I&\lorw&W^{\sharp}_I)\\
((\exp X,Y)&,&(\vu,\vrho,\vtheta))&\longmapsto&(\exp X\cdot\vu,\vrho,\vtheta+Y)
\end{array}.
$$
Analogously for the $\Lambda_J$-action.
Remark that the projections $W^{\sharp}_I\lorw\tilde{W}_I$ and 
$W^{\sharp}_J\lorw\tilde{W}_J$ induce homeomorphisms $W^{\sharp}_I/\Lambda_I\cong W_I$ and
$W^{\sharp}_J/\Lambda_J\cong W_J$.
We now  exhibit an equivariant homeomorphism $g^{\sharp}_{IJ}$ that projects down to
$g_{IJ}$.
This is given by the following mapping:
$$
\begin{array}{cccc}
g^{\sharp}_{IJ}:&W^{\sharp}_I&\lorw&W^{\sharp}_J\\
&(\vu,\vrho,\vtheta)&\longmapsto&
\exp(\pi_J^{-1}\cdot\pi)(\vtheta)\cdot\left(\vu+
(F_I(\vu,\sqrt{\vrho}\exp\vtheta))_{_{(J\setminus I\cap J)}}\right)
\end{array}.
$$
It is straightforward to check that this is a 
continuous, injective mapping
between open subsets of $\C^n$ whose
Jacobian matrix has
rank $2n$ at every point, therefore $g^{\sharp}_{IJ}$ is a 
diffeomorphism. Now add all compatible charts to obtain a complete atlas.
The key point that allows us to construct the lift $g^{\sharp}_{IJ}$ is 
the following: we can choose, for a point in $U_I\cap U_J$,
two  representatives, one in the slice corresponding to the open neighborhood
$\uci$ and one in the slice corresponding to the open neighborhood $\ucj$.
The mapping $g^{\sharp}_{IJ}$ expresses how to go from the first representative
to the other by moving along the orbit  that joins
the two slices.

{\it Part II: Singular strata}.
Singular strata are easier to describe as quasifolds, since
each of them is covered by one chart.
But let first take care of vertices: for each singular vertex $\mu$ of $\D$,
the stratum ${\cal T}_{\mu}=(\zset\cap {\C^*}^{{\mu}^c})/N^{\mu}$ is a point,
since by Lemma~\ref{torustrick} the group $N^{\mu}$ acts transitively
on $\zset\cap {\C^*}^{{\mu}^c}$.
Now let $F$ be a singular face of $\D$ of dimension $p>0$, let 
 $\mu$ be a vertex contained in $F$ and let 
 $I\in{\cal I}_{\mu}$ such that
$\hbox{card}(I\cap I_F)=(n-p)$.
The quotient 
$$\zset \cap {\C^*}^{F^c}/N$$
is a quasifold covered by one chart,  defined in
the following way.
Consider the open subset of $(\C^*)^p$:
$$\begin{array}{rcl}\check{U}_{F,I}&=&\{\vw\in(\C^*)^{I\setminus(I\cap I_F)}
\cong(\C^*)^p
\;|
\sum_{h\in I\setminus(I\cap I_F)}a_{hl}|w_{h}|^2>0,\;\;\hbox{for}\;l\in I_{\mu}\setminus 
(I\cup I_F)\\
&&\sum_{I\cap I_F}a_{hr}\lambda_h+
\sum_{h\in I\setminus(I\cap I_F)}a_{hr}(|w_{h}|^2+\lambda_{h})-\lambda_r>0
\;\;\hbox{for}\;r\notin I_{\mu}\}.\end{array}$$
and the mapping 
\begin{equation}\label{modellolocalesing}
\begin{array}{ccc}
\check{U}_{F,I} &\lorw &
{\cal T}_F\\
\vw & \longmapsto & [\vw+\check{\phi}_{F,I}(\vw)]
\end{array}
\end{equation}
with
$$\check{\phi}_{F,I}(\vw)_k=0,\;\quad k\in I\cup I_F$$
$$\check{\phi}_{F,I}(\vw)_l=\sqrt{\sum_{h\in I\setminus (I\cap I_F)}
a_{hl}|w_{h}|^2},\;\quad l\in I_{\mu}\setminus 
(I\cup I_F)$$
and
$$\check{\phi}_{F,I}(\vw)_r=\sqrt{\sum_{h\in I\cap I_F}a_{hr}\lambda_h+
\sum_{h\in I\setminus (I\cap I_F)}
a_{hr}(|w_{h}|^2+\lambda_{h})-\lambda_r}\;
\quad r\notin I_{\mu}.$$
The quotient
$
\check{U}_{F,I}/\check{\G}_{I\setminus(I\cap I_F)}$
is a model, in the sense of Remark~\ref{simplyc}:
the induced mapping
from
$
\check{U}_{F,I}/\check{\G}_{I\setminus(I\cap I_F)}$
to
${\cal T}_F$
is continuous, it is  
injective since $N\cap T^{(I\cup I_F)}=\G_I N^F.$ It is also
surjective and closed: this
can be easily proved by applying the 
same arguments of Part I, Steps (I,b) and (I,c).
The open set $\check{U}_{F,I}$ is not simply connected,
therefore, as we have 
done for some of the charts
of the regular piece, we have to consider its universal covering 
in order to have a proper chart.
We then add all compatible charts in order to obtain a complete atlas. 
In particular the charts corresponding to those $J$'s in ${\cal I}$ 
satisfying the 
hypothesis of Lemma~\ref{torusspliteffe} are compatible.

{\it Part III: The decomposition.} 
It remains to check that the described pieces,
indexed by the singular faces of $\D$ plus the open face,
give a decomposition of $M$ according to 
Definition~\ref{decomposition}. 
It is easy to verify that all pieces are locally closed and connected.
 The set of indexes certainly satisfies
point (i) of Definition~\ref{decomposition}. 
Point (ii) follows from Proposition~\ref{phi}, while point
(iii) is a consequence of Remark~\ref{stratidizset}.
Moreover the regular piece is open since the set of regular points in
$\zset$ is open, it is also dense since it contains the dense set
$\Phi^{-1}(\hbox{Int}(\D))$. 
\qed
\begin{remark}{\rm Each point $m$ in the space $M$ lies in a stratum, 
therefore it 
follows from Remark~\ref{gruppino} that there is a well defined
discrete group $\G_m$ attached to it.}\end{remark}
\begin{remark}{\rm
A singular face has at most dimension $n-3$, therefore a singular
piece has at most dimension $2n-6$}.\end{remark}
\begin{remark}\label{strutturalisciaglobale}{\rm The decomposition
of $M$ is induced by the decomposition of $\zset$ given by the
manifolds $\zset\cap {\C^*}^{F^c}$, with $F$ singular, and the open
subset $\zset_{\D}$ of $\zset$. The quasifold structure of each
piece ${\cal T}_F$ is naturally induced by the smooth structure of
$\zset\cap {\C^*}^{F^c}$, the quasifold structure of ${\cal
T}_{\D}$ is induced by the smooth structure of
$\zset_{\D}$}.\end{remark} Let $p\,\colon\,\zset\lorw M$ be the
projection, we have the following
\begin{thm}[Symplectic structure of strata]\label{stratisimplettici} Each piece ${\cal T}_F$ (${\cal T}_{\D}$)
of the decomposition of $M$ has a natural symplectic structure
induced by the quotient procedure, that is, its pull-back via $p$
coincides with the restriction of the standard symplectic form of
$\C^d$ to the manifold $\zset\cap {\C^*}^{F^c}$ ($\zset_{\D}$).\end{thm}
\proof Consider the regular piece. 
As in the classical reduction procedure, the standard symplectic 
structure of $\C^d$ induces a symplectic structure on each slice,
and therefore, via pullback, a symplectic structure, $\G_I$ ($\Lambda_I$)-- 
invariant, on each open subset $\uti\subset\C^n$ ($U^{\sharp}_I\subset\C^n$). 
The structure induced is the standard one and respect the changes of charts, 
thus defining a symplectic structure on the quasifold ${\cal T}_{\D}$.
The proof for the singular pieces goes in the same way.
\qed
\begin{thm}[Quasitorus action on $M$]\label{momentmapping} The restriction of the $D$-action
and of the mapping $\Phi$ to each piece of the space $M$ is
smooth, the action of $D$ is Hamiltonian and a moment mapping is
given by the restriction of $\Phi$.\end{thm}
\proof We refer to \cite{p} for the definition of Hamiltonian action
of a quasitorus on a quasifold and for the definition of moment mapping
with respect to this action. To prove that the action $\tau$ defined in
(\ref{action}) is smooth
and Hamiltonian we have to prove that it is so when lifted to local models.
From Proposition~\ref{phi} it then follows that the restriction of $\Phi$ to 
each stratum is a moment mapping.
We consider first the regular stratum ${\cal T}_{\D}$.
For each $I\in{\cal I}$, we have that the following diagram commutes
$$
\begin{array}{ccc}
\d\times \uti & \stackrel{\tilde{\tau}_{I}}{\mbox{\LARGE $\longrightarrow$}} &
\uti \\ (X,\vu)&{\mbox{\LARGE
$\longmapsto$}} & \exp (\pi_{I}^{-1}(X))\cdot \vu\\
{}\!\!\!\!\!\stackrel{}{}\,\stackrel{}{\mbox{\LARGE $\downarrow$}}
&  & \stackrel{}{\mbox{\LARGE $\downarrow$}}\,\stackrel{}{}\\
(\d/Q)\times(\uti/\G_I) & \stackrel{\tau_{I}}{\mbox{\LARGE $\longrightarrow$}} &
\uti/\G_I\\ ([X],[\vu])&{\mbox{\LARGE $\longmapsto$}} & [\exp
(\pi_{I}^{-1}(X))\cdot \vu]\\{}\!\!\!\!\!\stackrel{}{}\,\stackrel{}{\mbox{\LARGE
$\downarrow$}} &  & \stackrel{}{\mbox{\LARGE $\downarrow$}}\,\stackrel{}{}\\
(\d/Q)\times U_I & \stackrel{\tau}{\mbox{\LARGE $\longrightarrow$}} &
U_I=\widehat{U}_I/N\\
([X],\phi_I([\vu]))&{\mbox{\LARGE $\longmapsto$}} & 
[\exp(\pi_{I}^{-1}(X))\cdot\vu+F_I(\vu)]
\end{array}
$$ moreover $\tilde{\tau}_I$ is a smooth mapping and the action is Hamiltonian
with respect to the standard symplectic form on  $\uti$.
For each singular piece ${\cal T}_F$ we proceed in the same manner, consider
the diagram 
$$
\begin{array}{ccc}
\d\times \check{U}_{F,I} & \stackrel{\tilde{\tau}_{F,I}}{\mbox{\LARGE $\longrightarrow$}} &
\check{U}_{F,I} \\ (X,\vw)&{\mbox{\LARGE
$\longmapsto$}} & \exp (\pi_{I}^{-1}(X))_1\cdot \vw\\
{}\!\!\!\!\!\stackrel{}{}\,\stackrel{}{\mbox{\LARGE $\downarrow$}}
&  & \stackrel{}{\mbox{\LARGE $\downarrow$}}\,\stackrel{}{}\\
(\d/Q)\times(\check{U}_{F,I}/\check{\G}_{I\setminus(I\cap I_F)}) & 
\stackrel{\tau_{F,I}}{\mbox{\LARGE $\longrightarrow$}} &\check{U}_{F,I}
/\check{\G}_{I\setminus(I\cap I_F)}\\ ([X],[\vw])&{\mbox{\LARGE $\longmapsto$}} & [\exp
(\pi_{I}^{-1}(X))_1\cdot \vw]\\{}\!\!\!\!\!\stackrel{}{}\,\stackrel{}{\mbox{\LARGE
$\downarrow$}} &  & \stackrel{}{\mbox{\LARGE $\downarrow$}}\,\stackrel{}{}\\
(\d/Q)\times {\cal T}_F & \stackrel{\tau}{\mbox{\LARGE $\longrightarrow$}} &
{\cal T}_F\\
([X],[\vw+\check{\phi}_{F,I}(\vw)])&{\mbox{\LARGE $\longmapsto$}} & [\exp
(\pi_{I}^{-1}(X))_1\cdot \vw+\check{\phi}_{F,I}(\vw))]
\end{array}
$$
where $\exp
(\pi_{I}^{-1}(X))_1$ stands for the first factor of $\exp
(\pi_{I}^{-1}(X))$ in the decomposition $T^I=T^{I\setminus I\cap I_F}\times
T^{I\cap I_F}$.
The diagram  is commutative, the mapping
 $\tilde{\tau}_{F,I}$ is smooth and the action is Hamiltonian
with respect to the standard symplectic form on  $\check{U}_{F,I}$.
\qed
\section{The stratification: local structure}
Now we need to prove that our decomposition has a good local
behavior. Let $t_0$ be a point in the singular $2p$-dimensional
piece ${\cal T}_F$; we want to construct a link of $t_0$ satisfying
Definition~\ref{stratificazione}. Let
$\C^F$ 
be the complexification of the Lie algebra
${\R}^{F}$. The mapping $\Upsilon$ restricted to
${\C}^{F}$ gives rise to a moment mapping $\Upsilon_F$ for
the action of $T^F$ on  ${\C}^{F}$. The mapping
$\Upsilon_F\,\colon\,{\C}^{F}\lorw
\left({\R}^{F}\right)^*$ is given by
$\Upsilon_F(\vz)=\sum_{j\in I_F}(|z_j|^2+\lambda_j)e^*_j$.
Consider now the
Hamiltonian action of the $(r_{F}-n+p)$-dimensional group $N^{F}$
on ${\C}^{F}$, induced by that of $T^F$: a
moment mapping is then given by $\Psi_F=\iota_F^*\circ \Upsilon_F$, where
$\iota_{F}\,\colon\,\n^{F}\lorw{\R}^{F}$ is the
inclusion map. Using (\ref{facce}) we find that
$\Psi_{F}(\vz)=\sum_{j\in I_F}|z_j|^2\iota_{F}^*(e^*_j)$, hence
$\Psi^{-1}_F(0)$ is a cone. 
Let $\jmath_F:\d_F\lorw\d$ be
the inclusion map.
Consider the exact sequence 
$$0\lorw(\d_F)^*\stackrel{\pi_F^*}{\lorw}(\R^F)^*
\stackrel{\iota_F^*}{\lorw}(\n^F)^*\lorw0$$ 
where $\pi_F$ is the restriction to $\R^F$ of the projection $\pi:\R^{d}\lorw\d$.
Define
$\Sigma^{\diamond}_F=\bigcap_{j\in I_F}\{\;\mu\in\ddu\;|\;
\langle\mu,X_j\rangle\geq\lambda_j\;\}$.
By repeating the argument
of Proposition~\ref{phi} a continuous surjective
mapping 
$\Phi_F:\Psi_F^{-1}(0)/N^F\lorw\jmath_F^*(\sigmaf)$ is defined by setting
$\Phi_F([\vz])=(\pi_F^*)^{-1}(\Upsilon_F(\vz))$.
\begin{remark}\label{coniconetti}{\rm 
We enumerate here properties of $\sigmaf$ that we 
need in the sequel.
\begin{enumerate}
\item[(1)] $\Sigma_F=\jmath_F^*(\Sigma^{\diamond}_F)$ is an 
$(n-p)$-dimensional cone with vertex $\jmath_F^*(F)$;
\item[(2)] if $G$ is a $q$-dimensional face of $\D$ containing $F$, 
then $\jmath_F^*(G)$ is a $(q-p)$-dimensional face of ${\Sigma_F}$;
\item[(3)] for each $j\in I_F$ we can find $b_j\in(0,1]$ such that,
taken $Y=\sum_{j\in I_F}b_j X_j$, the intersection
${\Sigma_F}\cap\{\xi\in\d_F^*\;|\;\langle\xi,Y\rangle=
\sum_{j\in I_F}\lambda_j b _j+\epsilon\}$ is a nonempty convex 
polytope, $\D_{F,\epsilon}$, of dimension $(n-p-1)$  ($\epsilon\in\R_{>0}$);
\item[(4)] let $G$ be a $q$-dimensional face of $\D$ 
properly containing $F$, then 
$\GF=\jmath_F^*(G)\cap\{\xi\in\d_F^*\;|\;\langle\xi,Y\rangle=
\sum_{j\in I_F}\lambda_j b _j+\epsilon\}$
is a $(q-p-1)$-dimensional face of $\D_{F,\epsilon}$. Moreover
$\GF$ is singular in $\D_{F,\epsilon}$ if and only if $G$ is singular
in $\D$.
\end{enumerate}}
\end{remark}
Fix 
$Y\in\d_F$ as specified in the 
previous remark,
then the following Lemma holds:
\begin{lemma}[Local structure]\label{lemmalink} For each $t_0$ in the singular 
piece
${\cal T}_F$, 
the space $L_{F,\epsilon}
=\Phi_F^{-1}(\D_{F,\epsilon})$ 
satisfies the first point of Definition~\ref{stratificazione}
for a suitable $\epsilon$.
\end{lemma}
\proof
Let ${\cal S}_{F,\epsilon}=
\{\vz\in\C^F\;|\;
\sum_{j\in I_F}b_j|z_j|^2=\epsilon\}$ and let
$$(\Psi_F^{-1}(0))_{\epsilon}=\Psi_F^{-1}(0)\cap\{\vz\in\C^F\;|\;
\sum_{j\in I_F}b_j|z_j|^2<\epsilon\}.$$ 
Notice that
$L_{F,\epsilon}$ is nothing but
the quotient
$\left(\Psi_F^{-1}(0)\cap {\cal S}_{F,\epsilon}\right)/N^F$
and therefore
$$(\Psi_{F}^{-1}(0))_{\epsilon}/N^F
=C(L_{F,\epsilon}).$$ The space $L_{F,\epsilon}$ is our candidate link
for $t_0$. Recall that 
the decomposition in pieces of the space $M$
reflects the geometry
of the polytope $\D$ and
is defined via the mapping
$\Phi$.
The decompositions in pieces of both
spaces $(\Psi^{-1}_{F}(0))_{\epsilon}/N^F$ and 
$L_{F,\epsilon}$ are defined via the mapping $\Phi_F$ exactly in the same way
and
they are related accordingly to 
Remark~\ref{cono}.
The arguments used in the 
proof of Theorems~\ref{stratiquasifold} and \ref{stratisimplettici} apply 
with no important changes
to show that the two decompositions satisfy
Definition~\ref{decomposition} and that
their pieces are 
quasifolds, 
symplectic
in the case of $(\Psi^{-1}_{F}(0))_{\epsilon}/N^F$. 
Consider for instance
the singular pieces of $L_{F,\epsilon}$.
Each singular piece corresponds to a 
singular face $G$ of $\D$ properly containing $F$. 
Let $q$ be
the dimension  of $G$.
We want to prove that 
the singular piece $\Phi_F^{-1}(\GF)$
is a quasifold 
of dimension
$(2q-2p-1)$, covered by one chart.
The polytope to be considered is now $\D_{F,\epsilon}$.
We can choose $I\in{\cal I}$ 
such that $\hbox{card}(I\cap I_G)=n-q$ and $\hbox{card}(I\cap I_F)=n-p$,
then, for the transversality condition of Remark~\ref{coniconetti}, point
(3), there exists a $j\in I\cap I_F\setminus I\cap I_G$ such that,
having set $c_h=b_h+
\sum_{k\in I_F\setminus(I_G\cup(I\cap I_F))}b_ka_{hk}$, 
the coefficient $c_j\neq0$. By
$\sum_{h\neq j}$ we shall mean $\sum_{h\in (I_F\cap I)\setminus((I_G\cap I)
\cup\{j\})}$.
Define the discrete groups $\check{\G}_{I\cap I_F\setminus I\cap I_G}$
and $\check{\G}^{I\cap I_F}_{I\cap I_G}$ in such a  way that
$\G_{I\cap I_F}=\check{\G}_{I\cap I_F\setminus I\cap I_G}\times
\check{\G}^{I\cap I_F}_{I\cap I_G}$.
As in Theorem~\ref{stratiquasifold}, in order to construct a local model
we have to define a suitable slice.
Consider
the open subset:
\begin{equation}\label{ucheck}
\begin{array}{rl}\check{U}^F_{G,I,j}=&\{
(\theta_j,\vw)\in\R\times(\C^*)^{[(I_F\cap I)\setminus(I_G\cap I)]
\setminus\{j\}}\;|\;
(1/c_j)(\epsilon-\sum_{h\neq j}c_h|w_h|^2)>0,\\
&\sum_{h\neq j}a_{hk}|w_h|^2+(a_{jk}/c_j)(\epsilon-\sum_{h\neq j}
c_h|w_h|^2)>0,\;\;\;k\in I_F\setminus(I_G\cup(I\cap I_F))\}\end{array}
\end{equation}
acted on by the group
$\Z\times\check{\G}_{I\cap I_F\setminus(I\cap I_G)}$
in the following way:
\begin{equation}\label{azionedizeta}
\begin{array}{ccc}
(\Z\times\check{\G}_{I\cap I_F\setminus(I\cap I_G)})
\times \check{U}^F_{G,I,j}&\lorw&\check{U}^F_{G,I,j}\\
\left((m,\exp(T)),(\theta_j,\vw)\right)&\longmapsto&
(\theta_j+m+T_j,\exp(T_{(I_F\cap I)\setminus\{j\}})\vw).
\end{array}
\end{equation}
The homeomorphism
from the model 
$\check{U}^F_{\G,I,j}/\Z\times\check{\G}_{I\cap I_F\setminus(I\cap I_G)}$
to $\Phi_F^{-1}(\GF)\subset L_{\GF,\epsilon}$
is induced by the continuous mapping:
\begin{equation}
\label{modellolocalesing2}
\begin{array}{ccc}
\check{U}_{\GF,I,j}&\lorw&\Phi_F^{-1}(\GF)\\
(\theta_j,\vw)&\longmapsto&[\vw+\check{\phi}_{G_F}(\theta_j,\vw)]\end{array}
\end{equation}
with
$$
\begin{array}{ll}
\check{\phi}_{G_F}(\theta_j,\vw)_l=0&l\in I_G\cup (I_F)^c\cup (I\cap I_F)\\
\check{\phi}_{G_F}(\theta_j,\vw)_j=\sqrt{
(1/c_j)
\left(\epsilon-
\sum_{h\neq j}
c_h
|w_h|^2\right)}
e^{(2\pi i \theta_j)}&\\
\check{\phi}_{G_F}(\theta_j,\vw)_k=\sqrt{\sum_{h\neq j}a_{hk}|w_h|^2+
(a_{jk}/c_j)(\epsilon-
\sum_{h\neq j}
c_h|w_h|^2)}&k\in I_F\setminus (I_G\cup(I\cap I_F).
\end{array}
$$
The proof that the mapping 
(\ref{modellolocalesing2}) induces 
a homeomophism onto $\Phi^{-1}_{F}(G_F)$
goes very similarly to that given for the mapping (\ref{modellolocalesing}),
in the proof of Theorem~\ref{stratiquasifold}, Part~I, 
but since we deal now with the group $N^F$, the 
key result here is Lemma~\ref{torusspliteffe}.
The atlas, obtained by adding all compatible charts, contains, in this 
case too, all of the charts corresponding to those $J\in{\cal I}$ and
$j\in J$ satisfying the conditions specified above.

{\it The mapping $h_F$}. Let $F$ be a singular face of dimension $p>0$ and 
let $t_0$ be a point in ${\cal T}_F$.
We prove that near to $t_0$ our space $M$ is homeomorphic to 
the twisted product of an
open subset of ${\cal T}_F$ by a cone over the link $L_{F,\epsilon}$.
Let $I\in{\cal I}$ such that
$\check{U}_{F,I}/\check{\G}_{I\setminus(I\cap I_F)}$ gives a model for
${\cal T}_F$
as constructed in the proof of Theorem~\ref{stratiquasifold}, Part~II. 
In what follows we identify ${\cal T}_F$ with this model.
The discrete group $\check{\G}_{I\setminus(I\cap I_F)}$ acts on
the quotient $(\Psi^{-1}_{F}(0))_{\epsilon}/N^F$ and on
the product $\check{U}_{F,I}\times((\Psi^{-1}_{F}(0))_{\epsilon}/N^F)$ 
in the following 
way:
$$\begin{array}{ccccc}
\check{\G}_{I\setminus(I\cap I_F)}&\times&
(\Psi^{-1}_{F}(0))_{\epsilon}/N^F&\lorw& 
(\Psi^{-1}_{F}(0))_{\epsilon}/N^F\\
(g&,&[\vz])&\longmapsto&[\ell_{I_F,I}(g)\vz]\end{array};
$$
$$\begin{array}{ccccc}
\check{\G}_{I\setminus(I\cap I_F)}&\times&\left(
\check{U}_{F,I}\times((\Psi^{-1}_{F}(0))_{\epsilon}/N^F)\right)&\lorw& 
\check{U}_{F,I}\times((\Psi^{-1}_{F}(0))_{\epsilon}/N^F)\\
(g&,&(\vw,[\vz]))&\longmapsto&(g\vw,\ell_{I_F,I}(g)[\vz])\end{array}
.$$
By making use of the explicit
atlases, 
it is straightforward to check
that the quotient
$\left(
\check{U}_{F,I}\times((\Psi^{-1}_{F}(0))_{\epsilon}/N^F)\right)/
\check{\G}_{I\setminus(I\cap I_F)}$ inherits the
decomposition  in strata of the product, and that each of these  strata 
has an induced  
quasifold structure. This holds for any 
open subset $\check{B}\subset\check{U}_{F,I}$ such 
that $\check{B}/\check{\G}_{I\setminus I\cap I_F}$ is a model
in the sense of Remark~\ref{simplyc}.
We want to choose now a $\epsilon$ and a 
$T^{I\setminus(I\cap I_F)}$-
invariant open subset $\check{B}$ of $\check{U}_{F,I}$ such that
$t_0\in \check{B}/\check{\G}_{I\setminus(I\cap I_F)}$
and the mapping $h_F$ from the twisted product
$$\left(\check{B}\times((\Psi^{-1}_{F}(0))_{\epsilon}/N^F)\right)/
\check{\G}_{I\setminus(I\cap I_F)}$$ to the open subset of $M$
$$\left(\zset\cap(\check{B}\times
\{\vz\in\C^F\;|\;\sum_{j\in I_F}b_j|z_j|^2<\epsilon\}\times
\C^{(I_F\cup I)^c})
\right)/N$$ given by
$$h_F(
[\vw,[\vz]])=
[\vw+\vz+(F_I(\vz_{(I\cap I_F)}+\vw))_{(I\cup I_F)^c}]
$$
is well defined and surjective.
Let $\vw_0\in\check{U}_{F,I}$ be a point that projects down to
$t_0$, namely $t_0=[\vw_0+\check{\phi}_{F,I}(\vw_0)]$. 
We can choose a positive real constant $c$ and 
a $\check{B}$ containing 
$\vw_0$ in such a way that
$(\check{\phi}_{F,I}(\vw)_l)^2>c$ for all $l\in 
(I\cup I_F)^c$ and for all $\vw\in\check{B}$.
Denote by $B$ the open neighborhood in ${\cal T}_F$
homeomorphic to $\check{B}/\check{\G}_{I\setminus(I\cap I_F)}$.
The open subset $B$ is a submodel in ${\cal T}_F$.
Choose now $\epsilon>0$ in such a  way that
for each $[\vz]\in(\Psi^{-1}_{F}(0))_{\epsilon}/N^F$ 
we have 
$$\sum_{h\in I\cap I_F}a_{hl}|z_{h}|^2>-c,\;\quad l\in  
(I\cup I_F)^c.$$
With these choices 
the mapping $h_F$ is well defined, continuous 
and injective. It is easy to check, via Lemma~\ref{torustrick},
the $h_F$ is surjective. Moreover a simple adaptation of the
argument used 
in
the proof of Theorem~\ref{stratiquasifold}, Step (I,c), shows that the 
mapping $h_F$ is closed.
Finally we observe that
by construction 
the mapping $h_F$ takes strata into strata and its restriction to
each stratum is a  diffeomorphism of quasifolds, according
to Definition~\ref{diffeo}.
For each
singular vertex $\mu$ the mapping $h_{\mu}$ is defined on the cone
$C(L_{\mu,\epsilon})$ and satisfies all the required properties provided
that
$\epsilon>0$ is chosen in such a way  that $\sum_{h\in I}a_{hr}|z_{h}|^2>-
(\sum_{h\in I}
a_{hr}\lambda_h-\lambda_r)$ for each $r\in I_{\mu}^c$.
\qed
\begin{lemma}[The link of the link]\label{locale}
Let $F$ be a singular face of the convex polytope $\Delta$ and $t_0$
be a point in ${\cal T}_F$. The compact space $L_{F,\epsilon}$ is
a link of $t_0$.
\end{lemma}
\proof 
We need to prove that the decomposition of the compact space
$L_{F,\epsilon}$, defined in Lemma~\ref{lemmalink}, 
is itself a stratification, i.e. it
satisfies the recursive 
Definition~\ref{stratificazione}. 
Let $$G_0\subset G_1\subset G_2\subset\cdots\subset G_k$$ be
a sequence of singular faces of $\D$ such that $\hbox{dim}(G_r)=q_r$, 
with $G_0=F$ 
and such that there are no singular faces containing $G_k$. 
Notice  that $I_{G_r}\subset I_{G_{r-1}}$ for each $r=1,\dots,k$.
Now take a sequence of $Y^r\in\d_{G_r}$ and 
a sequence of $\epsilon_r$ in order to obtain, accordingly to Remark~\ref{coniconetti}, 
a corresponding sequence of convex polytopes $\D_{G_r,\epsilon_r}$, 
all nonsimple but the last one.
Choose an $I\in{\cal I}$ and a sequence of indices $j_{r-1}\in 
(I\cap I_{G_{r-1}})\setminus (I\cap I_{G_r})$ such that
the following two conditions are verified: the first
is
$\hbox{card}(I\cap I_{G_r})=n-q_r$
for each $\,r=0,\dots,k$
; for $r=1,\dots,k$ set
$c^{r-1}_h=b^{r-1}_h+\sum_k
b^{r-1}_ka_{hk}$, where,
if $h\in (I\cap I_{G_{r-1}})\setminus
(I\cap I_{G_r})$, then
$k$ ranges in ${I_{G_{r-1}}\setminus(I_{G_r}\cup(I\cap I_{G_{r-1}}}))$, 
if $h\in I\cap I_{G_r}$ then $k\in I_{G_{r-1}}\setminus I\cap I_{G_{r-1}}$,
the second condition then is
that the coefficient
$c^{r-1}_{j_{r-1}}\neq0$.
Let $t_r$ be a point in the singular piece ${\cal L}_{G_{r}}$ of
$L_{G_{r-1},\epsilon_{r-1}}$ corresponding to the face $G_r$. The space  
$L_{G_{r},\epsilon_{r}}$,
defined as
in Lemma~\ref{lemmalink} for a point in the piece ${\cal T}_{G_r}$, is
the candidate link of $t_r$.
Notice that $L_{G_k,\epsilon_k}$ is a quasifold.
The proof of the theorem is complete if we can 
prove that, for each such point $t_r$
, with $r=1,\dots,k$, the space 
$L_{G_{r-1},\epsilon_{r-1}}$
satisfies the first point of Definition~\ref{stratificazione}, namely
we need to prove that, near to $t_r$, the space $L_{G_{r-1},\epsilon_{r-1}}$
is homeomorphic to the twisted product of an open
subset of ${\cal L}_{G_r}$ by a cone over the link $L_{G_r,\epsilon_r}$.  
In order to do so
define $\check{U}^{G_{r-1}}_{G_r,I,j_{r-1}}$ in analogy with (\ref{ucheck}).
Consider now the discrete groups
$\check{\G}_{I\cap I_{G_{r-1}}\setminus I\cap I_{G_r}}$
and $\check{\G}^{I\cap I_{G_{r-1}}}_{I\cap I_{G_r}}$ such that
$\G_{I\cap I_{G_{r-1}}}=\check{\G}_{I\cap I_{G_{r-1}}\setminus I\cap I_{G_r}}
\times \check{\G}^
{
I\cap I_{G_{r-1}}
}
_{
I\cap I_{G_r}
}$.
The group $\Z\times\check{\G}_{I\cap I_{G_{r-1}}\setminus I\cap I_{G_r}}$
acts on $\check{U}^{G_{r-1}}_{G_r,I,j}$ as indicated in (\ref{azionedizeta}),
it also acts on $(\Psi_{G_r}^{-1}(0))_{\epsilon_r}/N^{G_r}$ 
in the following way:
$$\begin{array}{rclcc}
(\Z\times\check{\G}_{I\cap I_{G_{r-1}}\setminus I\cap I_{G_r}})
&\times&(\Psi_{G_r}^{-1}(0))_{\epsilon_r}/N^{G_r}&\lorw&
(\Psi_{G_r}^{-1})_{\epsilon_r}(0)/N^{G_r}\\
((m,g)&,&[\vz])&\longmapsto&\ell(g)[\vz]
 \end{array}
$$
where $\ell$ is the natural epimorphism 
$\ell\colon\check{\G}_{I\cap I_{G_{r-1}}\setminus I\cap I_{G_r}}
\lorw \check{\G}^{I\cap I_{G_{r-1}}}_{I\cap I_{G_r}}/\G_{I\cap I_{G_r}}$.
As in the proof of Lemma~\ref{lemmalink}
we choose 
an open neighborhood 
$\check{B}_r$ in $\C^{[(I\cap I_{G_{r-1}})\setminus 
(I\cap I_{G_r})]\setminus\{j_{r-1}\}}$, 
invariant by the action of  $T^d$, 
such that $\R\times\check{B}_r$ is contained in 
$\check{U}^{G_{r-1}}_{G_r,I,j_{r-1}}$ and the quotient
$\R\times\check{B}_r/(\Z\times\check{\G}_{I\cap I_{G_{r-1}}\setminus 
I\cap I_{G_r}})$ contains $t_r$. 
Consider now
the
mapping $h_r$ from the twisted product
$$\left(\R\times\check{B}_{r}\times 
(\Psi_{G_r,\epsilon_r}^{-1}(0)/N^{G_r})
\right)/\Z\times\check{\G}_{I\cap I_{G_{r-1}}\setminus I\cap I_{G_r}}$$
to the open subset of $L_{G_{r-1},\epsilon_{r-1}}$
$$\left(
\Psi_{G_{r-1}}^{-1}(0)\cap{\cal S}_{G_{r-1},\epsilon_{r-1}}\cap
\left(\check{B}_r\times\C^{\{j_{r-1}\}\cup I_{G_r}\cup (I\cap I_{G_{r-1}})^c}\right)\right)
/N^{G_{r-1}},$$ given by
$$h_r([\theta_j,\vw,[\vz]])=[\vw+\vz+\vx]$$
where
$$
x_h=0\quad\quad\hbox{for all}\quad h\in
[(I_{G_r}\cup(I\cap I_{G_{r-1}}))\setminus\{j_{r-1}\}]\cup (I_{G_{r-1}})^c
,$$
$$
x_{j_{r-1}}=
\sqrt{
(1/c^{r-1}_{j_{r-1}})
(\epsilon_{r-1}-
\sum_{h\neq j_{r-1}}
c^{r-1}_h|(\vw+\vz)_h|^2)}
e^{(2\pi i \theta_j)}$$
and 
$$x_s=
\sqrt{\sum_{h\neq j_{r-1}}a_{hs}|(\vw+\vz)_h|^2+
(a_{js}/c^{r-1}_{j_{r-1}})(\epsilon_{r-1}-
\sum_{h\neq j_{r-1}}
c^r_h|(\vw+\vz)_h|^2)}$$
where $h\neq j_{r-1}$ ranges in $I_{G_{r-1}}\setminus (I\cap I_{G_{r-1}})$
and 
$s\in 
 I_{G_{r-1}}\setminus((I\cap I_{G_{r-1}})\cup I_{G_{r}})
.$
The sequence of neighborhoods $\check{B_r}$, $r=1,\dots,k$, 
and the sequence of 
$\epsilon_r>0$, $r=0,\dots,k$, can be chosen in such a way that 
the mapping $h_r$ is well defined for each $r=1,\dots,k$.
A straightforward adaptation of the arguments used to check the properties
of the mapping $h_F$ in Lemma~\ref{lemmalink} shows that $h_r$ is 
continuous, bijective and closed.
Moreover $h_r$, restricted to each stratum, is a quasifold 
diffeomorphism.
\qed
\begin{thm}\label{teorema}
The decomposition of the space $M$ is a stratification by quasifolds
according to Definition~\ref{stratificazione}.
\end{thm}
\proof The proof is an immediate consequence of Lemma~\ref{lemmalink}
and Lemma~\ref{locale}.
\begin{remark}\label{linkovertoric}{\rm
The link $L_{F,\epsilon}$ fibers naturally over a space corresponding
to the polytope $\D_{F,\epsilon}$, this can be proved as follows:
let $\hbox{ann}(Y)=\{\xi\in\d_F^*\;|\;\langle \xi,Y\rangle=0\}$
and let $k_F:\hbox{ann}(Y)\lorw\d_F^*$ be the natural inclusion. 
Fix a point $\xi_0\in\{\xi\in\d_F^*\;|\;\langle\xi,Y\rangle=
\sum_{j\in I_F}\lambda_j b _j+\epsilon\}$ and denote by 
$\D'_{F,\epsilon}$ the polytope $\D_{F,\epsilon}$ viewed in the subspace 
$\hbox{ann}(Y)$. We have
$\D'_{F,\epsilon}=
\bigcap_{j\in I_F}\{\;\xi\in\hbox{ann}(Y)\;|\;
\langle\xi,k_F^*(X_j)\rangle\geq\lambda_j-\langle\xi_0,X_j\rangle\;\}
$.
Now  apply
the construction described in Section~3 to
the polytope $\D'_{F,\epsilon}$, 
with the choice of normals $k_F^*(X_j)$ and quasilattice 
$k_F^*(\d_F\cap Q)$. 
Denote by
$(N^{F})_{Y}$, $(\n^{F})_{Y}$ and 
$\Psi_{F,Y,\epsilon}$ the group thus obtained, its Lie algebra
and the relative moment mapping respectively. Let
$\tilde{Y}=\sum_{j\in I_F}b_je_j$, then it is straightforward to check that:

(i) $(\n^{F})_{Y}=\n_F\oplus \hbox{Span}\{\tilde{Y}\}$

(ii) $(N^{F})_{Y}/N^F\cong\exp(\hbox{Span}\{\tilde{Y}\})$ 

(iii) the moment mapping, written in components according to
the direct sum (i), is given by 
$\Psi_{F,Y,\epsilon}(\vz)=(\Psi_F(\vz),\sum_{j\in I_F}b_j|z_j|^2-\epsilon)$.
Therefore the link $L_{F,\epsilon}$ is exactly the quotient
$\Psi_{F,Y,\epsilon}^{-1}(0)/N^F$, 
hence $L_{F,\epsilon}$ fibers  
over the symplectic quotient
$\Psi_{F,Y,\epsilon}^{-1}(0)/(N^F)_{Y}$, with fiber
the $1$-dimensional group $(N^{F})_{Y}/N^F$.
Although the coefficients $b_j$'s can be chosen to be rational, 
thus obtaining a compact
fiber, this is not always the most natural choice, as we shall see in the
examples.}\end{remark}
\begin{remark}\label{choices}{\rm 
If the polytope $\D$ is rational, namely if the
$X_j$'s can be chosen to be in a lattice $L$,
then all discrete groups involved become finite, and therefore,
as we have already observed in Remark~\ref{quandobanale},
the stratification becomes locally trivial, singular strata are smooth
and the principal stratum is either smooth or an orbifold.
Each link $L_{F,\epsilon}$ is in this case a fiber bundle, with fiber
a 1-dimensional torus, over the
symplectic stratified space corresponding to the polytope
$\D'_{F,\epsilon}$, as specified in Remark~\ref{linkovertoric}. 
If in addition the polytope $\D$ admits a choice 
of $X_j$'s and
$Q$ that satisfies conditions (i) and (ii) specified in 
Remark~\ref{ifhonestlattice},
then the principal  stratum is also smooth.
The  quotient 
$\zset/N$ provides, in the rational case, an explicit
example of the symplectic stratified spaces described in \cite{ls},
as far as we refine the stratification of the principal stratum,
considering the strata given by finite group isotropy type. 
However  
the results obtained by Sjamaar-Lerman in \cite{ls} do 
not seem extendable to our specific case,
since they are based on Mather's results, which, as they are known,
do not apply to our context.

If the polytope $\D$ is simple each corresponding space has no singular strata
and  we find exactly the family of symplectic
quasifolds constructed by Prato in \cite{p}. This family includes, 
when $\D$ is rational in a lattice $L$, the symplectic orbifolds
associated with the pair $(\D,L)$, constructed in \cite{lt}; 
it also includes, when $\D$ satisfies Delzant's integrality condition, 
the symplectic toric manifold constructed in \cite{delzant}.
The simple case is significant in that it makes already clear that
orbifold structures are naturally associated to rational polytopes,
this explains why the principal stratum, which is associated with 
the regular part of the polytope, is, in general, an orbifold. 
As in the simple case, additional conditions
have to be satisfied in order to have a smooth principal stratum.
}\end{remark}
We are now ready to work out in detail Examples~\ref{piramide} and 
\ref{tenda}. For detailed examples of quasifolds, and in particular,
of the symplectic quasifolds corresponding to simple convex polytopes,
we refer the reader to \cite{p}.

\medskip

{\bf Example~\ref{piramide} resumed}.  
The regular stratum is a symplectic quasifold of dimension $6$.
It is covered by the open sets $U_{\mu_j}$, for $j=1,\dots,4$.
The corresponding models are $\tilde{U}_{\mu_j}/\G_{\mu_j}$.
The only singular 
stratum is the point ${\cal T}_{\nu}=[0,0,0,0,z_5]$. 
We want to describe $M$ as a cone
in a neighborhood of ${\cal T}_{\nu}$.
Let $I=\{2,3,4\}\subset{\cal I}_{\nu}$. The corresponding matrix is 
$$A_I=\left(\begin{array}{ccccc}
1/p_2&1&0&0&-p_5/p_2\\
-1&0&1&0&0\\
1&0&0&1&-p_5\end{array}\right).
$$ 
The subset $\zset\subset\C^5$ 
is given by the following equations:
$$|z_1|^2-1/p_2|z_2|^2-|z_4|^2+|z_3|^2=0,$$
$$|z_5|^2+p_5/p_2|z_2|^2+p_5|z_4|^2-p_5=0$$
while the $2$-dimensional group $N=\exp(\n)$ is the following subgroup of 
$T^5$:
$$N=\{(e^{2\pi i x},e^{2\pi i (-(1/p_2)x+(p_5/p_2)y)},e^{2\pi i x},
e^{-2\pi i(x+p_5y)},e^{2\pi i y})\;|\;x,y\in\R\}.$$
We construct now the link of ${\cal T}_{\nu}$:
the cone $\Psi_{\nu}^{-1}(0)/N^{\nu}$ is the quotient of
$$\{(z_1,z_2,z_3,z_4,0)\;|\;|z_1|^2-1/p_2|z_2|^2-|z_4|^2+|z_3|^2=0\}$$ 
modulo the action of
the $1$-dimensional group $$\begin{array}{r}N^{\nu}=N\cap T^{{\nu}}=\{(
e^{2\pi ix},e^{2\pi i(-1/p_2x+p_5/p_2h)},e^{2\pi i x},e^{2\pi i(-x+p_5 h
)},1)\;|\\
h\in\Z,x\in\R\}.\end{array}$$
Choose $\tilde{Y}=(1,1,1,p_2)$ according to 
Remarks~\ref{coniconetti}, ~\ref{linkovertoric}
and define the mapping 
$$\begin{array}{cccc}
h_{\nu}\colon&\Psi_{\nu,\epsilon}^{-1}(0)/N^{\nu}&\lorw&M\\
&[\vz]&\longmapsto&[z_1,z_2,z_3,z_4,\sqrt{p_5(1-|z_2|^2-|z_4|^2)}\,]
\end{array}
$$
where $\epsilon$ is chosen in such a way that if  
$\sum_{j=1}^{3}|z_j|^2+p_2|z_4|^2<\epsilon$, then $$-|z_2|^2-|z_4|^2>-1.$$
The link $L_{\nu,\epsilon}=(\Psi^{-1}_{\nu}(0)\cap {\cal S}_{\nu,\epsilon})
/N^{\nu}$ is a compact quasifold: it is the quotient
$$\{(z_1,z_2,z_3,z_4)\;|\;(1+1/p_2)|z_2|^2+2|z_4|^2=\epsilon\,;\;
|z_1|^2+|z_3|^2=\epsilon/(1+p_2)\}/N^{\nu}.$$
The group $N^{\nu}$ is a subgroup of the $2$-dimensional group 
$(N^{\nu})_{Y}=N_1\times N_2$ where $$N_1=\{\{(
e^{2\pi ix},1,e^{2\pi i x},1
,1)\;|\\
x\in\R\}$$ and 
$$N_2=\{(
1,e^{2\pi i(-1/p_2y+p_5/p_2h)},1,e^{2\pi i(-y+p_5 h
)},1)\;|\\
h\in\Z,y\in\R\}.$$ 
The quotient $(N^{\nu})_{Y}/N^{\nu}$ is isomorphic to
$\exp(\hbox{Span}\{\tilde{Y}\})$.
Therefore the link  fibers over $S^3/N_1\times S^3/N_2$, which is a 
product of two 
quasispheres (cf. \cite{p} for exact definitions and details). 
The fiber is the $1$-dimensional group
$(N^{\nu})_{Y}/N^{\nu}$, which is nonclosed if $p_2$ is nonrational.
Recall from \cite{p} that  
quasispheres 
are 
symplectic quasifolds associated to an interval in $\R$, therefore the link
$L_{\nu,\epsilon}$ consistently corresponds to $\D_{\nu,\epsilon}$, which, with
our choice of  $\tilde{Y}$, 
is a product of intervals.
If the $p_j$'s are rational we obtain a cone over an orbifold. If the
$p_j$'s are all equal to $1$, then conditions (i) and (ii) of 
Remark~\ref{ifhonestlattice} are satisfied and 
the corresponding space is stratified
by smooth manifolds. 
Near to the singular
point, the space is
a cone over  $L_{\nu,\epsilon}$, which is in this case
a fiber bundle over $S^2\times S^2$ with fiber $S^1$.

\medskip

{\bf Example~\ref{tenda} resumed}.
The regular stratum is an $8$-dimensional quasifold, a collection of charts
is given by the open subsets $U_I$, with $I$ ranging in ${\cal I}$.
We want to describe $M$ in a neighborhood of the singular 
point ${\cal T}_{\nu_1}=
[0,0,0,0,z_5,0,0,z_8,z_9]$ and in a neighborhood of a point in
the $2$-dimensional singular stratum ${\cal T}_{\nu_1\nu_2}$.
Consider $I=\{1,2,3,6\}\subset{\cal I}_{\nu_1}$. 
The corresponding matrix is 
$$A_I=\left(\begin{array}{ccccccccc}
1&0&0&1/p_1&0&0&-1/p_1&0&-1/p_1\\
0&1&0&1&-p_5&0&0&0&-1\\
0&0&1&-1&0&0&1&-p_8&0\\
0&0&0&0&-p_5&1&1&-p_8&0\end{array}\right).
$$ 
The subset $\zset\subset\C^9$ 
is given by the following equations:
$$|z_4|^2-1/p_1|z_1|^2-|z_2|^2+|z_3|^2=0,\quad\quad 
|z_7|^2+1/p_1|z_1|^2-|z_3|^2-|z_6|^2=0,$$
$$|z_5|^2+p_5|z_2|^2+p_5|z_6|^2-p_5=0,\quad\quad
|z_8|^2+p_8|z_3|^2+p_8|z_6|^2-p_8=0,$$
$$|z_9|^2+1/p_1|z_1|^2+|z_2|^2-1=0,$$
while the $5$-dimensional group $N=\exp(\n)$ is the following subgroup of 
$T^9$:
$$\begin{array}{r}N=\exp\{(1/p_1(-x+z+w),-x+p_5y+w,
x-z+p_8t,x,y,p_5y-z+p_8t,
z,t,w)\;|\;\\
x,y,z,t,w\in\R\}.\end{array}$$
We construct now the link of ${\cal T}_{\nu_1}$:
the cone $\Psi_{\nu_1}^{-1}(0)/N^{\nu_1}$ is the quotient of
$$\begin{array}{r}\{(z_1,z_2,z_3,z_4,0,z_6,z_7,0,0)\;|\;
|z_4|^2-1/p_1|z_1|^2-|z_2|^2+|z_3|^2=0,\\
|z_7|^2+1/p_1|z_1|^2-|z_3|^2-|z_6|^2=0\}\end{array}
$$ 
modulo the action of
the $2$-dimensional group $N^{\nu_1}=N\cap T^{{\nu_1}}$ given by
$$\begin{array}{r}\exp\{(1/p_1(-x+z+l),-x+p_5h+l,
x-z+p_8k,x,h,p_5h-z+p_8k,z,k,l)\;|\;\\
x,z\in\R\,\;h,k,l\in\Z\}.
\end{array}
$$
We choose $\tilde{Y}=(1,1,1,1,1,1)$ according to Remarks~\ref{coniconetti},
~\ref{linkovertoric}
and define the mapping 
$$\begin{array}{cccc}
h_{\nu_1}\colon&\Psi_{\nu_1,\epsilon}^{-1}(0)/N^{\nu_1}&\lorw&M\\
&[\vz]&\longmapsto&[\underline{z}']
\end{array}
$$
where 
$$\begin{array}{lcr}
z'_j=z_j&\hbox{for}&j\in I_{\nu_1},\\[.1cm]
z'_5=\sqrt{-p_5|z_2|^2-p_5|z_6|^2+p_5},&&\\[.1cm]
z'_8=\sqrt{-p_8|z_3|^2-p_8|z_6|^2+p_8},&&\\[.1cm]
z'_9=\sqrt{-1/p_1|z_1|^2-|z_2|^2+1}&&\end{array}$$
and $\epsilon$ is chosen in such a way that if  
$\sum_{j\in I_{\nu_1}}|z_j|^2<\epsilon$ then 
$$-|z_2|^2-|z_6|^2>-1,\quad\quad
-|z_3|^2-|z_6|^2>-1,$$
$$-1/p_1|z_1|^2-|z_2|^2>-1.$$
The link $L_{\nu_1,\epsilon}=
(\Psi^{-1}_{\nu_1}(0)\cap {\cal S}_{\nu_1,\epsilon})
/N^{\nu_1}$ is a $7$-dimensional stratified space. 
The singular pieces are $1$-dimensional
and correspond to the singular polytope edges stemming from $\nu_1$, namely:
$\nu_1\nu_2,\nu_1\mu_1$ and $\nu_1\mu_4$. They are the quotients
$$\{(0,0,0,0,0,z_6,z_7,0,0)\;|\;
|z_7|^2=|z_6|^2=\epsilon/2\}/N^{\nu_1},
$$ 
$$\{(0,z_2,0,z_4,0,0,0,0,0)\;|\;
|z_4|^2=|z_2|^2=\epsilon/2\}/N^{\nu_1},
$$ 
$$\{(z_1,0,z_3,0,0,0,0,0,0)\;|\;
1/p_1|z_1|^2=|z_3|^2=\epsilon/(1+p_1)\}/N^{\nu_1}.
$$ 
The mapping $h_{\nu_1}$ maps strata into strata 
diffeomorphically.

Let us now consider the $2$-dimensional singular stratum 
${\cal T}_{\nu_1\nu_2}$. Recall that $I_{\nu_1\nu_2}=\{1,2,3,4\}$ and that 
we have chosen
$I=\{1,2,3,6\}$. Following the proof of Theorem~\ref{stratiquasifold}
we can construct
a local model in a neighborhood of a point $t_0$ of the stratum: it is given by
$\check{U}_{I,\nu_1\nu_2}=\{w_6\in\C^*\;|\;|w_6|^2<1\}$
modulo the free action of the discrete group 
$\check{\G}_{I\setminus I\cap I_{\nu_1\nu_2}}$; denote by
$w_0$ the point in $\check{U}_{I,\nu_1\nu_2}$ projecting down to $t_0$.
The group $\check{\G}_{I\setminus I\cap I_{\nu_1\nu_2}}$
is obtained
by considering
$$\begin{array}{r}\G_I=\{(
e^{2\pi i1/p_1(-h+l+r)},e^{2\pi i(-h+p_5k+r)},
e^{2\pi i(h-l+p_8m)},1,1,\\e^{2\pi i(p_5k-l+p_8m)},
1,1,1)\;|\;
h,k,l,m,r\in\Z\}\end{array}$$
then 
$$\check{\G}_{I\setminus I_{\nu_1\nu_2}\cap I}=
\{(
1,1,1
,1,1,e^{2\pi i(p_5k-l+p_8m)},
1,1,1)\;|\;
k,l,m\in\Z\}.
$$
We also have the discrete group
$$\begin{array}{r}\check{\G}^I_{I\cap I_{\nu_1\nu_2}}=\{
e^{2\pi i1/p_1(-h+l+r)},e^{2\pi i(-h+p_5k+r)},
e^{2\pi i(h-l+p_8m)},1,1,1,
1,1,1)\;|\;\\
h,k,l,m,r\in\Z\}.\end{array}$$
Moreover, if $p_5,p_8/p_5\in\R\setminus\Q$ then 
$$\begin{array}{r}
{\G}_{I\cap I_{\nu_1\nu_2}}=\{
(e^{2\pi i(1/p_1)(-h+l+r)},1,
1,1,1,1,
1,1,1)\;|\;
h,l,r\in\Z\}.\end{array}$$
Recall that a natural group epimorphism is defined from
the group $\check{\G}_{I\setminus I_{\nu_1\nu_2}\cap I}$
onto the group 
$\check{\G}^I_{I\cap I_{\nu_1\nu_2}}/\G_{I\cap I_{\nu_1\nu_2}}$.
The cone $\Psi_{\nu_1\nu_2}^{-1}(0)/N^{\nu_1\nu_2}$ is the quotient of
$$\{(z_1,z_2,z_3,z_4,0,0,0,0,0)\;|\;
|z_4|^2-1/p_1|z_1|^2-|z_2|^2+|z_3|^2=0\}
$$ 
by the action of
the $1$-dimensional group $N^{\nu_1\nu_2}=N\cap T^{{\nu_1\nu_2}}$ given by
$$\begin{array}{r}\{(e^{2\pi i1/p_1(-x+l+r)},e^{2\pi i(-x+r)},
e^{2\pi i(x-l)},e^{2\pi i x},1,e^{-2\pi il},e^{2\pi il},1,e^{2\pi ir}\;|\;\\
x\in\R\,\;l,r\in\Z\}\end{array}$$
We choose $\tilde{Y}=(1,p_1,1,1)$ according to Remarks~\ref{coniconetti},
~\ref{linkovertoric} 
and define the mapping 
$$\begin{array}{cccc}
h_{\nu_1\nu_2}\colon&(\tilde{B}\times\Psi_{\nu_1\nu_2,\epsilon}^{-1}(0)/
N^{\nu_1\nu_2})\,/\,\check{\G}_{I\setminus(I\cap I_{\nu_1\nu_2})}&\lorw&M\\
&[\vw,[\vz]]&\longmapsto&[\underline{z}']
\end{array}
$$
where
$$\begin{array}{lcr}
z'_j=z_j&\hbox{for}&j\in I_{\nu_1\nu_2}\\
z'_6=w_6,&&\\[.1cm]
z'_5=\sqrt{-p_5|z_2|^2-p_5|w_6|^2+p_5},&&\\[.1cm]
z'_7=\sqrt{-1/p_1|z_1|^2+|z_3|^2+|w_6|^2},&&\\[.1cm]
z'_8=\sqrt{-p_8|z_3|^2-p_8|w_6|^2+p_8},&&\\[.1cm]
z'_9=\sqrt{-1/p_1|z_1|^2-|z_2|^2+1}
\end{array}$$ 
and a positive constant $c$ is chosen in such a way that
$\tilde{B}=\{w_6\in\C^*\;|\;c<|w_6|^2<1-c
\}$ is well defined and contains $w_0$,
$\epsilon$ is chosen in such a way that if  
$|z_1|^2+p_1|z_2|^2+|z_3|^2+|z_4|^2<\epsilon$ then 
$$-|z_2|^2>-c,\quad\quad 
-|z_3|^2>-c,\quad\quad
-1/p_1|z_1|^2-|z_2|^2>-c,$$
$$-1/p_1|z_1|^2+|z_3|^2>-c.$$
Here we can touch the twisting group
$\check{\G}_{I\setminus I\cap I_{\nu_1\nu_2}}$. 
Recall from Lemma~\ref{torustrick}
that $N/\G_I\cong T^d/T^I$; namely $\G_I$, when infinite, 
represents, intuitively, the 
nonclosed part of $N$.  The twisting group
$\check{\G}_{I\setminus I\cap I_{\nu_1\nu_2}}$ 
is a subgroup of $\G_I$ and it does
act on both sides of the product 
$(\tilde{B}\times\Psi_{\nu_1\nu_2,\epsilon}^{-1}(0)/N^{\nu_1\nu_2})$.
The link $L_{\nu_1\nu_2,\epsilon}=
(\Psi^{-1}_{\nu_1\nu_2}(0)\cap {\cal S}_{\nu_1\nu_2,\epsilon})
/N^{\nu_1\nu_2}$ is a compact quasifold which fibers  
over the product
of two quasispheres, with fiber a $1$-dimensional 
group isomorphic to $\exp(\hbox{Span}\{\tilde{Y}\})$,
similarly to the link found in the pyramid example.

To exemplify the proof of Theorem~\ref{locale} let us
consider the sequence of singular faces $\nu_1\subset\nu_1\nu_2$.
The corresponding polytopes, $\D_{\nu_1,\epsilon_1}$ and
$\D_{\nu_1\nu_2,\epsilon_2}$, 
can be visualized in Figure~3, 
\begin{figure}[h]
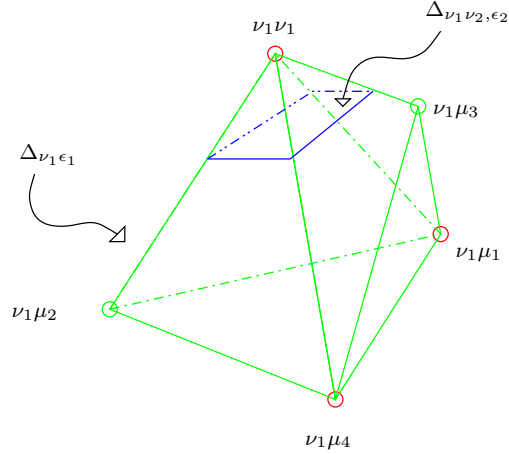

\begin{center}
\input polylinks.pstex_t
\end{center}
\caption{The link of the link}
\end{figure}
the link of 
the singular point ${\cal T}_{\nu_1}$, namely
$L_{\nu_1,\epsilon_1}$, is a fibration over a space corresponding
to $\D_{\nu_1,\epsilon_1}$,  
the link of the link, at a singular point
in the stratum corresponding to $\nu_1\nu_2$, is a fibration over
a space corresponding to $\D_{\nu_1\nu_2,\epsilon_2}$.
The fibers are, as we have seen, $1$-dimensional abelian group, possibly 
nonclosed.
If the $p_j$'s are rational the space corresponding to our polytope is 
stratified by 
orbifolds. If the
$p_j$'s are all equal to $1$, then conditions (i) and (ii) of 
Remark~\ref{ifhonestlattice} are satisfied and the corresponding space is, 
near
to each singular stratum, a trivial bundle over the stratum itself; moreover
the strata are smooth  manifolds.
\begin{remark}\label{strutture}{\rm
Theorem~\ref{teorema} proves that the decomposition of
$M$ is in fact a stratification. 
Moreover, from
Theorems~\ref{stratiquasifold},$\,$\ref{stratisimplettici}, we know that
each piece of the stratification of $M$ has the structure of a
symplectic quasifold, naturally induced by that of $\C^d$.
This suggests 
that the symplectic forms defined on each stratum glue together to 
give rise to a 
symplectic form on the stratified space $M$, globally defined, thus 
making sense of the notion of differential form defined on $M$.
}
\end{remark}
\begin{remark}{\rm In the light of Theorem~\ref{momentmapping},
we can view the mapping $\Phi$ as a moment mapping for the action of
the $n$-dimensional
quasi-torus $D$ on the $2n$-dimensional compact space $M$ stratified
by symplectic quasifolds. 
By Proposition~\ref{phi}, the image $\Phi(M)$ of the moment
mapping $\Phi$, is exactly
the polytope $\D$.}\end{remark}
\begin{remark}\label{complex}{\rm The remark above
emphasizes the relationship between the space $M$ and the polytope
$\D$, which is very neat in the symplectic setting. From the
complex point of view we have a compact space $X$, 
stratified by complex quasifolds and a homeomorphism from $M$ onto $X$
that is a diffeomorphism restricted to the strata. The space
$X$ is $n$-dimensional and
is acted on by the complexified torus $\Dc$ of same dimension.
Such an action has a dense open orbit, corresponding to the open
set $\Phi^{-1}(\hbox{Int}(\D))$. Complex toric spaces corresponding to $\D$ 
will be treated in \cite{nscx}.}\end{remark}

\noindent \small{\sc Dipartimento di Matematica Applicata ``G. Sansone'',
Via S. Marta 3, 50139 Firenze, ITALY,  
{\tt mailto:fiammetta.battaglia@unifi.it}}
\end{document}